\documentclass[12pt]{article}
\usepackage{graphics}
 \usepackage{tikz}
 \usepackage{booktabs}
 \usetikzlibrary{calc}
\usepackage{graphicx}
\usepgflibrary{bbox}
\usepackage{pgfplots} 
\usepackage{mathrsfs,pstricks}
\usepackage{amsmath,amsfonts,amsthm,amssymb}
\usepackage{epic,eepic,epsfig}
\usepackage{latexsym, enumerate}
\usepackage{graphics,multicol}
\usepackage{color} 
 \newtheorem{thm}{Theorem}[section]
 \newtheorem{cor}[thm]{Corollary}

 \theoremstyle{definition}

 \newtheorem{defn}[thm]{Definition}

\newcommand{\proofend}[1]{\hfill$\diamondsuit$ #1\smallskip}

 \renewcommand\tilde{\widetilde}

\begin{document}
\baselineskip 0.22 in

\title{Nonorientable genus embedding of nearly complete bipartite graphs\thanks{Supported by the NNSFC, China under Grant No. 11301171 and by Hunan Province Natural Science Fund Projects, China under
Grant No. 2019JJ40080.}}
\author{Shengxiang Lv\thanks{E-mail address: lvsxx23@126.com} 
 \\{\footnotesize School of Big Data and Statistics,
Hunan University of Finance and Economics,}\\ {\footnotesize  Changsha 410205, P.R.China}\\ 
 } 
\date{}
\maketitle
 
 \noindent{\bf Abstract}\quad  The nearly complete bipartite graph $G(m,n,k)$ is obtained by removing  $k$ independent edges from the complete bipartite graph $K_{m,n}$. In this paper, we prove that  for any nearly complete bipartite graph $G(m,n,k)$ with $m, n\geq 3$,  and $(m,n,k)\notin\{(5,4,4)$, $(4,5,4)$, $(5,5,5)\}$, there exists  a nonorientable genus embedding $\Pi$ satisfying $\tilde{\gamma}(\Pi)=\max\{\lceil \big((m-2)(n-2)-k\big)/2\rceil, 1\}$. This embedding can be constructed by starting from an embedding of some $G(p,q,h)$ with $h\leq 6$ and $p,q\leq 7$, and then iteratively adding multiple copies of $G(2,2,2)$, $G(2,0,0)$ and $G(0,2,0)$.  As a consequence, the previously unresolved nonorientable genus $\tilde{\gamma}(G(n+1,n,n))$ for even $n$ and $\tilde{\gamma}(G(n,n,n))$ for arbitrary $n$ are now determined.
 
  \bigskip
 

 
 \vskip 0.1cm

\section {Introduction}
Let $G$ be a graph with vertex set $V(G)$ and edge
set $E(G)$. The complete bipartite graph with partition sets $A$ and $B$ is denoted by $K_{A,B}$. When $|A|=m$ and $|B|=n$, we also write $K_{A, B}$ as $K_{m,n}$. For convenience, we use $G_{A,B}$ to denote
a bipartite subgraph of the complete biparttite graph $K_{A,B}$.  
The \emph{nearly complete bipartite graph} $G(m,n,k)$ is obtained
 from the complete bipartite graph $K_{m,n}$ by deleting $k$ independent edges, where $k\leq \min\{m,n\}$. 
Clearly, $G(m,n,0)$ is the complete bipartite graph $K_{m,n}$. 
    
A surface is a connected compact $2$-manifold without boundary. 
Let $S_p$ denote the \emph{orientable surface} with \emph{genus} $p$, 
 and $N_q$ the \emph{nonorientable surface} with \emph{nonorientable genus} $q$. The \emph{Euler characteristic} of a surface $\Sigma$, denoted $\chi (\Sigma)$, is $2-2p$ for $S_p$, and $2-q$ for $N_q$. 
The \emph{Euler genus} of  $\Sigma$ is defined as $\mathcal{E}(\Sigma)=2-\chi (\Sigma)$. 
 
An \emph{embedding} of a graph $G$ on a surface  $\Sigma$ is a drawing of $G$ on $\Sigma$ such that no edge passes through any vertex and no two edges do not cross. An embedding $\Pi$ of a graph on a surface $\Sigma$ is called \emph{cellular} if every component of $\Sigma-\Pi$ is homeomorphic to an open disc. The \emph{Euler genus} of such an embedding $\Pi$ is defined as $\mathcal{E}(\Pi)=\mathcal{E}(\Sigma)$. Throughout this paper, all embeddings considered are cellular.
  
  In a cellular embedding $\Pi$, each component of $\Sigma-\Pi$ is called
a \emph{face}, and the closed walk bounding a face is called a \emph{facial walk}.  If the closed walk $a_{1}a_2 \cdots a_{k}a_1$ is the  facial walk of a face $f$, then we write $f=(a_{1}a_2 \cdots a_{k})$.
   
For an embedding $\Pi$ on a surface $\Sigma$, if $\Sigma$ is orientable, then $\Pi$ is called an \emph{orientable embedding}, and its \emph{orientable genus}  is defined as $\gamma(\Pi)=\frac{1}{2}\mathcal{E}(\Sigma)$. If $\Sigma$ is nonorientable, then $\Pi$ is called a \emph{nonorientable embedding}, and its \emph{nonorientable genus} is $\widetilde{\gamma}{(\Pi)}=\mathcal{E}(\Sigma)$.
   
The  \emph{orientable genus} of a graph $G$, denoted by $\gamma(G)$, is the smallest 
 integer $p\geq 0$ such that $G$ can be embedded in the orientable surface $S_{p}$. Any such embedding is called an \emph{orientable genus embedding} of $G$.
 
 The \emph{nonorientable genus} of $G$, denoted by $\widetilde{\gamma}(G)$, is the smallest 
 integer $q>0$ such that $G$  can be embedded in the nonorientable surface $N_{q}$.
Any such embedding is called a \emph{nonorientable genus embedding} of $G$.
 
 An embedding in which every face is bounded by a $4$-cycle is called a \emph{quadrilateral embedding} 
 (or \emph{quadrangular}), also known as a \emph{quadrangulation} of the surface. If all faces except one
  satisfy this condition, the embedding is called a \emph{nearly quadrilateral embedding}. Similarly, a 
  \emph{triangular embedding} (or \emph{triangulation}) is an embedding in which every face is bounded by a $3$-cycle. 
  
   The neighborhood $N(v)$ of a vertex $v$ of a graph $G$ is the set of vertices adjacent to $v$.  A cyclic permutation is written as $[x_1,x_2,\cdots,x_t]$.
    Since we restrict our attention to simple graphs, a \emph{rotation system}  for $G$ can be represented by a set $\{\rho_v: v\in V(G)\},$ where $\rho_v$ is a cyclic permutation of  $N(v)$, called the  \emph{rotation at~$v$}. If $D=\{\rho_v: v\in V(G)\}$ is a rotation system  of $G$, we denote $\rho_v$ by $D_v$ .
Two rotation systems $D$ and $D'$ of $G$ are \emph{equivalent} if $D'=\{D_v, v\in V(G)\}$ or $D'=\{(D_v)^{-1}, v\in V(G)\}$.
 
   Consider an embedding $\Pi$ of graph $G$.  An orientation at a vertex $v$ in $\Pi$ \emph{induces a rotation} $\Pi_v$.
  By  arbitrarily fixing an orientation at every vertex of  $\Pi$, the set of  induced rotations $\{\Pi_v, v\in V(G)\}$ defines a rotation system 
of $G$,  called the rotation system of $G$ \emph{induced} by $\Pi$. 

It is well known that nonequivalent rotation systems correspond to distinct cellular embeddings of the graph \cite{St1}.

 Determining the orientable or nonorientable genus of a graph is NP-complete\cite{Th1, Th2}.
A triangular embedding of a graph is a minimum genus embedding. Ringel, 
in his proof of the Map Color Theorem\cite{Ri1}, 
constructed a triangular embedding of the complete graph $K_n$. 

Quadrilateral embeddings have also been extensively studied. 
For  complete graphs, such embeddings were studied in \cite{Liu1, Z1, HT1, HT2}. 
Hartsfield and Ringel \cite{HT1, HT2} also discovered quadrilateral embeddings of the general octahedral
graph $Q_{2n}$. Hartsfield \cite{H1} provided a nonorientable quadrilateral
embedding of a complete multipartite graph $K_{n_1,n_2,\cdots,n_t}$(other than $K_5$ and $K_{1,m,n}$) with an even number of edges. Moreover, Quadrilateral embeddings for graphs obtained from certain graph operations were studied in \cite{D1, P1, P2, W1,W2}.

For bipartite graphs, minimum genus embeddings are often quadrilateral or nearly quadrilateral. For convenience, define
 $$ f(m,n,k)=\frac{(m-2)(n-2)-k}{2}.$$
 By Euler's formula, one easily obtains the following lower bound on the nonorientable genus 
\begin{equation}
\label{eq-1}
\tilde{\gamma}(G(m, n, k))\geq \lceil f(m,n,k)\rceil . 
\end{equation}
Equality holds if and only if the embedding is quadrilateral or nearly quadrilateral. 

Ringel\cite{Ri2, Ri3} and Bouchet\cite{Bo1} determined
the minimum genus of complete bipartite graphs.
 Mohar et al.\cite{M1, MPP1} studied the orientable and nonorientable genus  
of nearly complete bipartite by constructing quadrilateral  or nearly quadrilateral embeddings.   
\ In \cite{MPP1},  Mohar et al. determined the orientable genus of all nearly complete bipartite graphs $G(m, n, k)$.  
 The nonorientable genus of those graphs has not yet been fully determined.  For the following cases in Theorem \ref{N-1},  Mohar determined  that equality holds in Equation \eqref{eq-1}.
    
  {\thm [{Mohar \cite[Theorem 1]{M1}}]\label{N-1}  Let $m, n\geq 3$ and $k\leq \min\{m, n\}$. If the triple  $(m, n, k)$ is not one of the following : $(n + 1, n, n)$ and $(n, n + 1, n)$
   with $n$ even, or $(n, n, n)$ with $n$ arbitrary, then
$$\tilde{\gamma}[G(m, n, k)] = \max\{ \lceil  f(m, n, k)\rceil, 1\}.$$}   
   
For $G(5,4,4)$ and $G(5,5,5)$, Mohar \cite{M1}  further established that the inequality in Equation \eqref{eq-1} is strict, i.e.,
   $$\tilde{\gamma}[G(5, 4, 4)] =2~~~\textrm{and}~~~\tilde{\gamma}[G(5, 5, 5)] =3.$$
   Therefore, the cases $(n+1,n,n)$, $(n,n+1,n)$ with $n$ even, and $(n,n,n)$ for arbitrary $n$ remain open.
 
 Let $G_{A,B}$ and $G_{C,D}$ be two bipartite graphs such that $(A \cup B)\cap (C\cup D)=\emptyset$.  
The \emph{bipartite join} of $G_{A,B}$ and $G_{C,D}$, denoted by $G_{A,B}\oplus G_{C,D}$, is defined as  the bipartite graph with vertex set
 $$V=A\cup B\cup C\cup D$$
 and edge set 
 $$E=E(G_{A,B})\cup E(G_{C,D})\cup E(C,B)\cup E(A, D),$$
 where $ E(X,Y)$ denotes the set of all edges between vertex sets $X$ and $Y$.
 Equivalently, $G_{A,B}\oplus G_{C,D}$ is the bipartite graph with partition sets $A\cup B$ and $C\cup D$ and the edge set described above. 
 
 Clearly, the nearly complete bipartite graphs $$G(m+2,n,k),\ G(m,n+2,k),\ G(m+2,n+2,k+2)$$ can be obtained via the following bipartite join constructions:
 \begin{align*}
G(m+2,n,k)&=G(m,n,k)\oplus G(2,0,0),\\
 G(m,n+2,k)&=G(m,n,k)\oplus G(0,2,0),\\
 G(m+2,n+2,k+2)&=G(m,n,k)\oplus G(2,2,2).
\end{align*}
For notational simplicity, exponents will be used hereafter to denote repeated bipartite joins. For example, $G(m,n,k)=G(p,q,h) \oplus G(r,s,t)^x$ denotes 
$$G(m,n,k)=G(p,q,h) \oplus G(r,s,t) \oplus \cdots \oplus G(r,s,t)$$
 with $x$ copies of $G(r,s,t)$. Moreover, an extend embedding $\Pi(m,n,k)$ of $G(m,n,k)$ is introduced in Section \ref{sec3} (see Definition \ref{defn3-5})  and denoted as 
 $$\Pi(m,n,k)=\Pi(p,q,h)\oplus G(r,s,t)^x,$$
 where $\Pi(p,q,h)$ is an embedding of $G(p,q,h)$ satisfying $\mathcal{E}(\Pi(p,q,h))=\lceil f(p,q,h)\rceil$.

 \vskip0.3cm
   \begin{table}[htbp]\label{table1}
     \begin{center}
     \caption{The graphs $G(p,q,h)$ of Theorem \ref{thm2}.}
  \renewcommand{\arraystretch}{1.1}
  \begin{tabular}{ccccccc}
 \toprule[1pt]
  $h$       &$h=0$	        &$h=1$	       &$h=2$	     &$h=3$        &$h=4$        &$h=5,6$\\
\midrule[1pt]
               & $G(2,2,0)$	&$G(2,2,1)$     &$G(3,3,2)$  &$G(3,4,3)$   &$G(4,4,4)$  &$G(5,5,5)$\\
$G(p,q,h)$&$G(2,3,0)$	&$G(2,3,1)$	&$G(3,4,2)$ &$G(3,5,3)$   &$G(4,7,4)$	&$G(6,7,6)$\\
              &$G(3,3,0)$	&$G(3,3,1)$	& 	           &	                &                  &                \\
\bottomrule[1pt]
  \end{tabular}
\end{center}
  \label{table1}
  \end{table}
   
Employing an inductive procedure, we begin with an embedding $\Pi(p,q,h)$ of $G(p,q,h)$ with small $p$, 
 $q$ and $h$ (see Table \ref{table1}). By iteratively adding two or four new vertices, we obtain an extend embedding $\Pi(m,n,k)$ (see Definition \ref{defn3-5} in Section \ref{sec3}) of a larger nearly complete bipartite graph $G(m,n,k)$ satisfying $\tilde{\gamma}(\Pi(m,n,k))= \max\{ \lceil  f(m, n, k)\rceil, 1\}$. The larger nearly complete bipartite graph $G(m,n,k)$ is precisely the bipartite join of $G(p,q,h)$ with multiple copies of $G(2,0,0)$, $G(0,2,0)$, or $G(2,2,2)$. That is 
$$G(m,n,k)=G(p,q,h)\oplus G(2,2,2)^a\oplus G(2,0,0)^b\oplus G(0,2,0)^c,$$
where 
$$p+2a+2b=m,\ q+2a+2c=n, \ h+2a=k.$$ The first main result of this paper is given by the following theorem.

  {\thm\label{thm2} Let $G(m,n,k)$ be a nearly complete bipartite graph with $m,n\geq 3$, and $(m,n,k)\notin\{(5,4,4)$, $(4,5,4)$, $(5,5,5)\}$, then there exist a  nearly complete bipartite graph $G(p,q,h)$ from Table \ref{table1} such that 
  $$G(m,n,k)=G(p,q,h)\oplus G(2,2,2)^a\oplus G(2,0,0)^b\oplus G(0,2,0)^c,$$
  where 
$$p+2a+2b=m,\ q+2a+2c=n, \ h+2a=k.$$
 Moreover,  $G(p,q,h)$ admits an embedding  $\Pi(p,q,h)$ with $\mathcal{E}(\Pi(p,q,h))= \lceil f(p, q, h)\rceil $ such that the extend embedding
  $$\Pi(m,n,k)=\Pi(p,q,h)\oplus G(2,2,2)^a\oplus G(2,0,0)^b\oplus G(0,2,0)^c$$
  is a nonorientable embedding of $G(m,n,k)$ satisfying
  $$\tilde{\gamma}(\Pi(m,n,k))= \max\{ \lceil  f(m, n, k)\rceil, 1\}.$$
  } 
 
The nonorientable genus of $G(m,n,k)$ for all cases except $(m,n,k)\neq (5,4,4)$, $(4,5,4)$, $(5,5,5)$ follows directly from Theorem \ref{thm2} and Equation \ref{eq-1}. Consequently, the previously open cases on $\tilde{\gamma}(G(n+1,n,n))$ for even $n$ and $\tilde{\gamma}(G(n,n,n))$ for arbitrary $n$ are now resolved.
  
  {\cor\label{N-2} For $m,n\geq 3$ with $(m,n,k)\neq (5,4,4)$, $(4,5,4)$, $(5,5,5)$,  the nonorientable genus of $G(m,n,k)$ is 
$$\tilde{\gamma}(G(m, n, k)) = \max\{\lceil  f(m, n, k)\rceil,1\}.$$
  In particular, for $n\geq 6$, 
$$\tilde{\gamma}(G(n, n, n)) = \lceil  f(n, n, n)\rceil=  \frac{n^2-5n+4}{2}.$$
and for even $n\geq 6$,  
$$\tilde{\gamma}(G(n+1, n, n))= \lceil  f(n+1, n, n)\rceil= \frac{n^2-4n+2}{2}.$$
}

This paper is organized as follows: In Section 2, we introduce the concept of an extendible pair in an embedding and extendible embedding of $G(m,n,k)$, then  prove that every embedding of the bipartite graph $G(m,n,k)$ are extendible.
In Section 3, we construct extended embeddings of $G(m,n,k)\oplus G(2,2,2)$, $G(m,n,k)\oplus G(2,0,0)$ and $G(m,n,k)\oplus G(0,2,0)$ from a given embedding of $G(m,n,k)$.  Finally, in Section 4, we present the proof of Theorem \ref{thm2}. 
 
For further background on graph genus and embeddings, we refer the reader to \cite{Gr1, MT1}.
 
\section{Extendible embedding of $G(m,n,k)$}\label{sec2}

Let $G_{A,B}$ be a connected bipartite graph with partition sets 
$$A=\{a_1,a_2,\ldots, a_m\} \text{ (red vertices),} \ B=\{b_1,b_2,\ldots,b_n\} 
\text{(black vertices)},$$
where $m,n\geq2$.

 Let $\Pi$ be a cellular embedding of $G_{A,B}$ on a surface $\Sigma$ and let $\mathcal{F}$ denote the set of faces of $\Pi$. If there exits a face $f\in \mathcal{F}$ of the form 
 $$f=(x_i \cdots x_j\cdots y_s\cdots y_t\cdots),$$ 
 where $x_i,x_j,y_s,y_t$ are four distinct vertices in $A\cup B$, then we say that the two pairs $\{x_i,x_j\}$ and $\{y_s,y_t\}$ are \emph{parallel in $f$}.  

 \begin {defn}  Given an embedding $\Pi$ of a bipartite graph $G_{A,B}$, let
$$\mathcal{F}_A=(f_1,f_2,\ldots, f_{\lceil m/2\rceil}) \text{ and } \mathcal{F}_B=(g_1,g_2,\ldots,g_{\lceil n/2\rceil}),$$
be two ordered sequences of faces in $\Pi$. If the faces in $\mathcal{F}_A$ and $\mathcal{F}_B$ satisfying the following form  (with vertices suitably relabeled if necessary):
  \begin{align*}
f_{\lceil m/2\rceil}&=g_{\lceil n/2\rceil}=(a_{m-1}\cdots b_{n-1} \cdots a_{m}\cdots  b_{n} \cdots),\\
f_i&=(a_{2i-1}\cdots b'_i a_{2i}\cdots b''_i),\quad  \text{for}\ 1\leq i\leq \lfloor (m-2)/2\rfloor, \\
g_j&=(b_{2j-1}\cdots a'_j b_{2j}\cdots a''_j), \quad  \text{for}\  1\leq j\leq  \lfloor (n-2)/2\rfloor, \\
\intertext{and, when $m$ is odd, }
 f_{\lceil (m-2)/2\rceil}= f_{(m-1)/2}&=(a_{2\lfloor (m-2)/2\rfloor}b'_{\lceil (m-2)/2\rceil}\cdots a_{2\lfloor(m-2)/2\rfloor+1}b''_{\lceil (m-2)/2\rceil}\cdots) \\
 &=(a_{m-3}b'_{(m-1)/2}\cdots a_{m-2}b''_{(m-1)/2}\cdots),
\intertext{when $n$ is odd, }
 g_{\lceil (n-2)/2\rceil}= g_{(n-1)/2}&=(b_{2\lfloor (n-2)/2\rfloor}a'_{\lceil (n-2)/2\rceil}\cdots b_{2\lfloor(n-2)/2\rfloor+1}a''_{\lceil (n-2)/2\rceil}\cdots) \\
 &=(b_{n-3}a'_{(n-1)/2}\cdots b_{n-2}a''_{\lceil(n-1)/2\rceil}\cdots ).
\end{align*}
where all $a'_j, a''_j \in A$ and $b'_i,b''_i \in B$,
then the pair $(\mathcal{F}_A,  \mathcal{F}_B)$ is called
 a \emph{saturated pair of face sequences}.
 \end{defn}
 
 In a saturated pair of face sequences, the set of the two named $a$ vertices in the face $f_i$  is called the \emph{$f_i$-pair}; similarly, 
the set of the two named $b$ vertices in the face $g_j$ is called the \emph{$g_j$-pair}. Thus,
 for $1\leq i\leq \lfloor (m-2)/2\rfloor$, the $f_i$-pair is $\{a_{2i-1}, a_{2i}\}$; 
 the $f_{\lceil m/2\rceil}$-pair is $\{a_{m}, a_{m-1}\}$; when $m$ is odd, the $ f_{({m-1})/2}$-pair is $\{a_{m-3}, a_{m-2}\}$.

In a saturated pair of face sequences $(\mathcal{F}_A,  \mathcal{F}_B)$, faces may be repeated in the ordered sequence $\mathcal{F}_A$ or $\mathcal{F}_B$, and a face may appear in both sequences. Although two faces in $(\mathcal{F}_A, \mathcal{F}_B)$ may be identical, their associated pairs are distinct.  For example, even though $f_{\lceil m/2\rceil}=g_{\lceil n/2\rceil}$, the $g_{\lceil n/2\rceil}$-pair is different from the $f_{\lceil m/2\rceil}$-pair. 

Moreover, when $m$ is odd, the vertex $a_{m-3}$ belongs to both the $f_{(m-3)/2}$-pair and the $f_{(m-1)/2}$-pair. Similarly, when $n$ is odd, the vertex $b_{n-3}$   belongs to both the $g_{(n-3)/2}$-pair and the $g_{(n-1)/2}$-pair.
 
 {\defn A saturated pair of face sequences $(\mathcal{F}_A, \mathcal{F}_B)$ is called an \emph{extendible pair} if the following four conditions are satsfied:
\begin{enumerate}
\item If $f_{s}=f_t$ for some $1\leq s< t\leq \lceil m/2\rceil$, then the $f_s$-pair and the $f_t$-pair are parallel in $f_{s}=f_t$.
\item If $g_{p}=g_q$ for $1\leq p< q\leq  \lceil n/2\rceil$, then the $g_p$-pair and the $g_q$-pair are parallel in $g_{p}=g_q$.  
\item  
If $f_i=g_j$ for some $1\leq i \leq  \lceil (m-2)/2\rceil$ and  $1\leq j \leq  \lceil n/2\rceil$, then the $f_i$-pair and  $g_j$-pair are parallel in $f_i=g_j$.
\item  
If $g_j=f_{ \lceil m/2\rceil}$ for some $1\leq j \leq  \lceil (n-2)/2\rceil$, then the $g_j$-pair and  $f_{\lceil m/2\rceil}$-pair are parallel in $g_{j}=f_{\lceil m/2\rceil}$.
\end{enumerate}
}

 {\defn An embedding $\Pi$ of $G_{A,B}$ is called \emph{extendible} if it admits an extendible pair $(\mathcal{F}_A, \mathcal{F}_B)$.}
 
 The following theorem states that, apart from a few exceptional cases, every embedding of a connected nearly complete bipartite graph 
 $G({m,n,k})$ is extendible.
 
\begin{thm}\label{extendible1}
For $m,n\geq 2$, if  $(m,n,k)\notin \{(2,2,2),(3,3,3), (4,4,4), (5,5,5)\}$, then every cellular embedding $\Pi$ of the nearly complete bipartite graph $ G({m,n,k})$ is extendible.
 \end{thm}
 
{\bf Proof}\quad Let $G(m,n,k)$ be the graph obtained from the complete bipartite graph $K_{A,B}$ by removing the edge set
$$M_k=\{a_ib_i, 1\leq i\leq k\}.$$ Thus, 
$$E(G(m,n,k))=E(K_{m,n})\setminus M_k.$$ We now prove, by considering three cases depending on the value of $k$,  that there  exists an extendible pair  $(\mathcal{F}_A, \mathcal{F}_B)$ in the embedding $\Pi$. 

{\bf Case 1.}\ $k< \min\{m,n\}$. Since $b_na_m\in E(G_{m,n,k})$, there is a face $f$  in $\Pi$ whose boundary contains $b_na_m$. Without loss of generality,  let $f=(a_{t_{m-1}}b_na_mb_{h_{n-1}}\cdots)$ and suppose the rotations at $b_n$ and $a_m$ induced by $\Pi$ are 
$$\Pi_{b_n}=[a_{t_1}a_{t_2}\cdots a_{t_{m-1}}a_{m}], ~~\Pi_{a_{m}}=[b_{h_1}b_{h_2}\cdots b_{h_{n-2}}b_nb_{h_{n-1}}],$$
as illustrated in Figure \ref{F1}.
 
 \begin{figure}[h]

\setlength{\unitlength}{1.5mm}
  \begin{center}
\begin{tikzpicture}[]

\fill[red!20,red!20] (1,0)--(1,1.5) -- (5,1.5)--(5,0)--(1,0);

 \coordinate [label={[label distance=0mm]-15:{\small$b_n$}}] (bn) at (1,1.5);
  \coordinate [label={[label distance=0mm]-135:{\small$a_m$}}]  (am) at (5,1.5);
\coordinate [label=below:{\small$a_{t_{m-1}}$}]  (am-1) at (1,0);
\coordinate [label=below:{\small $b_{n-1}$}]  (bn-1) at (5,0); 
 \fill[blue!20] (1,1.5)--({1+1.4*cos(30)}, {1.5+1.4*sin(30)})--({1+1.4*cos(60)}, {1.5+1.4*sin(60)})--(1,1.5);
\fill[blue!20] (1,1.5)--({1+1.4*cos(90)}, {1.5+1.4*sin(90)})--({1+1.4*cos(120)}, {1.5+1.4*sin(120)})--(1,1.5);
\fill[blue!20] (1,1.5)--({1+1.4*cos(170)}, {1.5+1.4*sin(170)})--({1+1.4*cos(200)}, {1.5+1.4*sin(200)})--(1,1.5);

\fill[red!20] (5,1.5)--({5+1.4*cos(-60)}, {1.5+1.4*sin(-60)})--({5+1.4*cos(-30)}, {1.5+1.4*sin(-30)})--(5,1.5);
\fill[red!20] (5,1.5)--({5+1.4*cos(0)}, {1.5+1.4*sin(0)})--({5+1.4*cos(30)}, {1.5+1.4*sin(30)})--(5,1.5);
\fill[red!20] (5,1.5)--({5+1.4*cos(70)}, {1.5+1.4*sin(70)})--({5+1.4*cos(100)}, {1.5+1.4*sin(100)})--(5,1.5);
 
\def\a{0};
\def\b{30};
\foreach \i in {1,2,3,4} {
\coordinate [label={[label distance=0mm]\a+\b*\i:{\small $a_{t_{\i}}$}}](a\i) at ({1+1.5*cos(\a+\b*\i)}, {1.5+1.5*sin(\a+\b*\i)});
\draw[thick] (bn) --(a\i) ;
\filldraw [draw=black, fill=red](a\i) circle (3pt);
}

\coordinate [label={[label distance=0mm]170:{\small $a_{t_{2i-1}}$}}](a5) at ({1+1.5*cos(170)}, {1.5+1.5*sin(170)});
\draw[thick] (bn) --(a5) ;
\coordinate [label={[label distance=0mm]200:{\small $a_{t_{2i}}$}}](a6) at ({1+1.5*cos(200)}, {1.5+1.5*sin(200)});
\draw[thick] (bn) --(a6) ;

\coordinate [label={[label distance=0mm]240:{\small $a_{t_{m-2}}$}}](a7) at ({1+1.5*cos(240)}, {1.5+1.5*sin(240)});
\draw[thick] (bn) --(a7) ;

\def\c{-90};
\def\d{30};
\foreach \i in {1,2,3,4} {
\coordinate [label={[label distance=0mm]\c+\d*\i:{\small $b_{h_{\i}}$}}](b\i) at ({5+1.5*cos(\c+\d*\i)}, {1.5+1.5*sin(\c+\d*\i)});
\draw[thick] (am) --(b\i) ;
\filldraw [draw=black, fill=blue] (b\i)circle (3pt);
}
\coordinate [label={[label distance=0mm]70:{\small $b_{h_{2j-1}}$}}](b5) at ({5+1.5*cos(70)}, {1.5+1.5*sin(70)});
\draw[thick] (am) --(b5) ;
\coordinate [label={[label distance=0mm]100:{\small $b_{h_{2j}}$}}](b6) at ({5+1.5*cos(100)}, {1.5+1.5*sin(100)});
\draw[thick] (am) --(b6) ;
\coordinate [label={[label distance=0mm]90:{\small $b_{h_{n-2}}$}}](b7) at ({5+1.5*cos(150)}, {1.5+1.5*sin(150)});
\draw[thick] (am) --(b7) ;

\filldraw [draw=black, fill=blue] (b5)circle (3pt)(b6)circle (3pt)(b7)circle (3pt);

\node at ({1+1*cos(45)}, {1.5+1*sin(45)}){$f_1$};
\node at ({1+1*cos(105)}, {1.5+1*sin(105)}){$f_2$};
\node at ({1+1*cos(185)}, {1.5+1*sin(185)}){$f_i$};

\node at ({5+1*cos(-45)}, {1.5+1*sin(-45)}){$g_1$};
\node at ({5+1*cos(15)}, {1.5+1*sin(15)}){$g_2$};
\node at ({5+1*cos(85)}, {1.5+1*sin(85)}){$g_j$};

\node at ({5+1*cos(40)}, {1.5+1*sin(40)}){$\cdot$};
\node at ({5+1*cos(50)}, {1.5+1*sin(50)}){$\cdot$};
\node at ({5+1*cos(60)}, {1.5+1*sin(60)}){$\cdot$};

\node at ({5+1*cos(115)}, {1.5+1*sin(115)}){$\cdot$};
\node at ({5+1*cos(125)}, {1.5+1*sin(125)}){$\cdot$};
\node at ({5+1*cos(135)}, {1.5+1*sin(135)}){$\cdot$};

\node at ({1+1*cos(135)}, {1.5+1*sin(135)}){$\cdot$};
\node at ({1+1*cos(145)}, {1.5+1*sin(145)}){$\cdot$};
\node at ({1+1*cos(155)}, {1.5+1*sin(155)}){$\cdot$};

\node at ({1+1*cos(210)}, {1.5+1*sin(210)}){$\cdot$};
\node at ({1+1*cos(220)}, {1.5+1*sin(220)}){$\cdot$};
\node at ({1+1*cos(230)}, {1.5+1*sin(230)}){$\cdot$};

\draw[very thick] (am-1)--(bn) -- (am)--(bn-1);
\draw[very thick,dotted] (bn-1) -- (am-1);
\filldraw [draw=black, fill=blue] (bn) circle (3pt)(bn-1) circle (3pt);
\filldraw [draw=black, fill=red](am) circle (3pt) (am-1) circle (3pt)(a5) circle (3pt)(a6) circle (3pt)(a7) circle (3pt);
\node at (3,0.75) {$f_{ \lceil m/2 \rceil}=g_{ \lceil n/2 \rceil}$};
\end{tikzpicture}
\end{center}

\caption{{\small The induced rotations at $b_n$ and $a_m$ for $k< \min\{m,n\}$.}}\label{F1}
 \end{figure}
 
We now choose the faces as follows: 
\begin{align*}
f_{ \lceil m/2 \rceil}=g_{ \lceil n/2 \rceil}=f&=(a_{t_{m-1}}b_na_mb_{h_{n-1}}\cdots),\\
f_{i}&=(a_{t_{2i-1}}b_na_{t_{2i}}\cdots)~\text{for}~1\leq i\leq \lfloor (m-2)/2 \rfloor, \\
g_{j}&=(b_{h_{2j-1}}a_mb_{h_{2j}}\cdots)~\text{for}~1\leq j\leq \lfloor (n-2)/2 \rfloor,\\
\intertext{and, if $m$ is odd, }
f_{ (m-1)/2}&=(a_{t_{m-3}}b_na_{t_{m-2}}\cdots), \\
\intertext{if $n$ is odd, }
g_{(n-1)/2}&=(b_{h_{n-3}}a_mb_{h_{n-2}}\cdots).
 \end{align*}
 
In particular, if $m=3$ (respectively, $n=3$), then 
$$f_1=(a_{t_{1}}b_na_{t_2}\cdots)\quad \big(\text{ resp.}~ g_1=(b_{h_{1}}a_mb_{h_2}\cdots)\big).$$  Consequently,
$$\big(\{f_1,f_2,\ldots, f_{\lceil m/2 \rceil}\},\{g_1,g_2,\ldots, g_{\lceil n/2 \rceil}\}\big)$$ 
is an extendible pair in $\Pi$.


\begin{figure}[h]

\setlength{\unitlength}{1.5mm}
  \begin{center}
\begin{tikzpicture}[]

\def\e {8.5};
\def\f{12.5};
\fill[red!20,red!20](\e,1.5)--(\f,1.5) -- (\f,3)--(\e,3)--(\e,1.5);
  \coordinate [label={[label distance=0mm]30:{\small$a_m$}}]  (am) at (\e,1.5);
   \coordinate [label={[label distance=0mm]95:{\small$b_{t_1}$}}] (bt1) at (\f,1.5);
   \coordinate [label=above:{\small $b_{t_2}$}]  (bt2) at (\e,3); 
\coordinate [label=above:{\small$a_{h_{m-1}}$}]  (ahm-1) at (\f,3);
\coordinate [label={[label distance=0mm]100:{\small$a_{h_x}=a_{t_2}$}}]  (at2) at (\f+4, 1.5);
 
\def\a{30};
\def\b{30};
\foreach \i in {3,4,5,6} {
\coordinate [label={[label distance=0mm]\a+\b*\i:{\small $b_{t_{\i}}$}}](b\i) at ({\e+1.5*cos(\a+\b*\i)}, {1.5+1.5*sin(\a+\b*\i)});
\draw[thick] (am)--(b\i) ;
\filldraw [draw=black, fill=blue](b\i) circle (3pt);
}

\coordinate [label={[label distance=0mm]260:{\small $b_{t_{n-2}}$}}](b7) at ({\e+1.5*cos(260)}, {1.5+1.5*sin(260)});
\draw[thick] (am)--(b7) ;
\coordinate [label={[label distance=0mm]300:{\small $b_{t_{n-1}}$}}](b8) at ({\e+1.5*cos(300)}, {1.5+1.5*sin(300)});
\draw[thick] (am)--(b8);
\filldraw [draw=black, fill=blue] (b7) circle (3pt)(b8) circle (3pt);

\node at({\e+1*cos(225)}, {1.5+1*sin(225)}){$\cdot$};
\node at ({\e+1*cos(235)}, {1.5+1*sin(235)}){$\cdot$};
\node at ({\e+1*cos(245)}, {1.5+1*sin(245)}){$\cdot$};
 
\def\a{180};
\def\b{30};
\foreach \i in {1,2,3,4} {
\coordinate [label={[label distance=0mm]\a+\b*\i:{\small $a_{h_{\i}}$}}](a\i) at ({\f+1.5*cos(\a+\b*\i)}, {1.5+1.5*sin(\a+\b*\i)});
\draw[thick] (bt1) --(a\i) ;
\filldraw [draw=black, fill=red](a\i) circle (3pt);
}

\coordinate [label={[label distance=0mm]30:{\small $a_{h_{2i-1}}$}}](a5) at ({\f+1.5*cos(30)}, {1.5+1.5*sin(30)});
\draw[thick] (bt1) --(a5) ;
\coordinate [label={[label distance=0mm]60:{\small $a_{h_{2i}}$}}](a6) at ({\f+1.5*cos(60)}, {1.5+1.5*sin(60)});
\draw[thick] (bt1) --(a6) ;

\fill[blue!20] (bt1)--({\f+1.4*cos(\a+\b*1)}, {1.5+1.4*sin(\a+\b*1)})--({\f+1.4*cos(\a+\b*2)}, {1.5+1.4*sin(\a+\b*2)})--(bt1);
\node at ({\f+1*cos(\a+45*1)}, {1.5+1*sin(\a+45*1)}){$f_1$};

\fill[blue!20] (bt1)--({\f+1.4*cos(\a+\b*3)}, {1.5+1.4*sin(\a+\b*3)})--({\f+1.4*cos(\a+\b*4)}, {1.5+1.4*sin(\a+\b*4)})--(bt1);
\node at ({\f+1*cos(\a+\b*3+15)}, {1.5+1*sin(\a+\b*3+15)}){$f_2$};

\fill[blue!20] (bt1)--({\f+1.4*cos(30)}, {1.5+1.4*sin(30)})--({\f+1.4*cos(60)}, {1.5+1.4*sin(60)})--(bt1);
\node at ({\f+1*cos(45)}, {1.5+1*sin(45)}){$f_i$};
\node at({\f+1*cos(320)}, {1.5+1*sin(320)}){$\cdot$};
\node at ({\f+1*cos(330)}, {1.5+1*sin(330)}){$\cdot$};
\node at ({\f+1*cos(340)}, {1.5+1*sin(340)}){$\cdot$};

\node at({\f+1*cos(10)}, {1.5+1*sin(10)}){$\cdot$};
\node at ({\f+1*cos(15)}, {1.5+1*sin(15)}){$\cdot$};
\node at ({\f+1*cos(20)}, {1.5+1*sin(20)}){$\cdot$};

\node at({\f+1*cos(70)}, {1.5+1*sin(70)}){$\cdot$};
\node at ({\f+1*cos(75)}, {1.5+1*sin(75)}){$\cdot$};
\node at ({\f+1*cos(80)}, {1.5+1*sin(80)}){$\cdot$};

\def\a{-150};
\def\b{30};
\foreach \i in {1,2,3,4} {
\coordinate [label={[label distance=0mm]\a+\b*\i:{\small $b_{l_{\i}}$}}](c\i) at ({\f+4+1.5*cos(\a+\b*\i)}, {1.5+1.5*sin(\a+\b*\i)});
\draw[thick] (at2) --(c\i) ;
\filldraw [draw=black, fill=blue](c\i) circle (3pt);
}
\coordinate [label={[label distance=0mm]10:{\small $b_{l_{2j-1}}$}}](c5) at ({\f+4+1.5*cos(10)}, {1.5+1.5*sin(10)});
\draw[thick] (at2) --(c5) ;
\coordinate [label={[label distance=0mm]30:{\small $b_{l_{2j}}$}}](c6) at ({\f+4+1.5*cos(40)}, {1.5+1.5*sin(40)});
\draw[thick] (at2) --(c6) ;
\coordinate [label={[label distance=0mm]80:{\small $b_{l_{n-2}}$}}](c7) at ({\f+4+1.5*cos(80)}, {1.5+1.5*sin(80)});
\draw[thick] (at2) --(c7) ;
\filldraw [draw=black, fill=blue](c5) circle (3pt)(c6) circle (3pt)(c7) circle (3pt);

\fill[red!20] (at2)--({\f+4+1.4*cos(\a+\b*1)}, {1.5+1.4*sin(\a+\b*1)})--({\f+4+1.4*cos(\a+\b*2)}, {1.5+1.4*sin(\a+\b*2)})--(at2);
\node at ({\f+4+1*cos(\a+\b*1+15)}, {1.5+1*sin(\a+\b*1+15)}){$g_1$};
\fill[red!20] (at2)--({\f+4+1.4*cos(\a+\b*3)}, {1.5+1.4*sin(\a+\b*3)})--({\f+4+1.4*cos(\a+\b*4)}, {1.5+1.4*sin(\a+\b*4)})--(at2);
\node at ({\f+4+1*cos(\a+\b*3+15)}, {1.5+1*sin(\a+\b*3+15)}){$g_2$};

\fill[red!20] (at2)--({\f+4+1.4*cos(10)}, {1.5+1.4*sin(10)})-- ({\f+4+1.4*cos(40)}, {1.5+1.4*sin(40)})--(at2);
\node at ({\f+4+1*cos(25)}, {1.5+1*sin(25)}){$g_j$};
\node at({\f+4+1*cos(-10)}, {1.5+1*sin(-10)}){$\cdot$};
\node at ({\f+4+1*cos(0)}, {1.5+1*sin(0)}){$\cdot$};
\node at ({\f+4+1*cos(-20)}, {1.5+1*sin(-20)}){$\cdot$};
\node at({\f+4+1*cos(50)}, {1.5+1*sin(50)}){$\cdot$};
\node at ({\f+4+1*cos(60)}, {1.5+1*sin(60)}){$\cdot$};
\node at ({\f+4+1*cos(70)}, {1.5+1*sin(70)}){$\cdot$};

\draw[very thick] (bt2)--(am)--(bt1) -- (ahm-1)(bt1) -- (at2);
\draw[very thick,dotted] (bt2) -- (ahm-1);
\filldraw [draw=black, fill=blue] (bt1) circle (3pt)(bt2) circle (3pt);
\filldraw [draw=black, fill=red](am) circle (3pt) (ahm-1) circle (3pt)(at2) circle (3pt)(a5) circle (3pt)(a6) circle (3pt);
\node at (10.5,2.4) {$f_{ \lceil m/2 \rceil}=g_{ \lceil n/2 \rceil}$};
\end{tikzpicture}
\end{center}

\caption{{\small The induced rotations at $a_m$, $b_{t_1}$ and $a_{t_2}$ for $2\leq k=m<n$, where $b_{t_1}a_{t_1}\notin M_m, b_{t_{2}}a_{t_{2}}\in M_m$.}}\label{F2}
 \end{figure}

{\bf Case 2.}\quad $2\leq k=m<n$. 
Suppose the rotation at vertex $a_m$ induced by $\Pi$ is
$$\Pi_{a_m}=[b_{t_1}b_{t_2}\cdots b_{t_{n-2}}b_{t_{n-1}}].$$
Because $k=m<n$, there exists an index $1\leq i\leq n-1$, without loss of generality, assume $i=1$, such that 
 $$b_{t_{1}}a_{t_{1}}\notin M_k,\ b_{t_{2}}a_{t_{2}}\in M_k.$$
Hence, $|N_{b_{t_1}}|=m$ and $|N_{a_{t_2}}|=n-1$,  and suppose the rotations at 
$b_{t_{1}}$ and $a_{t_2}$ induced by $\Pi$ are
$$\Pi_{b_{t_{1}}}=[a_{h_1}a_{h_2}\cdots a_{h_{m-2}}a_{h_{m-1}}a_{m}], \quad \Pi_{a_{t_{2}}}=[b_{l_1}b_{l_2}\cdots b_{l_{n-2}}b_{t_1}],$$
where $a_{h_x}=a_{t_{2}}$ for a $1\leq x\leq m-1$, and $$\{b_{l_1}, b_{l_2},\cdots, b_{l_{n-2}}, b_{t_1}\}=B\setminus \{b_{t_2}\}.$$
The rotations at $a_m$, $b_{t_{1}}$ and $a_{t_{2}}$ are illustrated in Figure \ref{F2}.
 
 Thus, we can choose the faces as follows: 
\begin{align*}
f_{ \lceil m/2 \rceil}=g_{ \lceil n/2 \rceil}&=(a_{h_{m-1}}b_{t_1}a_{m}b_{t_{2}}\cdots),\\
f_{i}&=(a_{h_{2i-1}}b_{t_1}a_{h_{2i}}\cdots) ~\text{for}~1\leq i\leq \lfloor (m-2)/2 \rfloor, \\
g_{j}&=(b_{l_{2j-1}}a_{t_2}b_{l_{2j}}\cdots) ~\text{for}~1\leq j\leq \lfloor (n-2)/2 \rfloor, \\
\intertext{and, if $m$ is odd, }
f_{(m-1)/2}&=(a_{h_{m-3}}b_{t_1}a_{h_{m-2}}\cdots),  \\
\intertext{if $n$ is odd, }
g_{(n-1)/2}&=(b_{l_{n-3}}a_{t_2}b_{l_{n-2}}\cdots).\\
 \end{align*}

In particular, when $2=k=m<n$, we have $a_m=a_2,$ $a_{t_2}=a_1$ and $b_{t_2}=b_1$. The chosen faces are
 \begin{align*}
f_{1}=g_{ \lceil n/2 \rceil}&=(a_{1}b_{t_1}a_{2}b_{1}\cdots),\\
g_{j}&=(b_{l_{2j-1}}a_{1}b_{l_{2j}}\cdots) ~\text{for}~1\leq j\leq \lfloor (n-2)/2 \rfloor, \\
\intertext{if $n$ is odd, }
g_{(n-1)/2}&=(b_{l_{n-3}}a_{1}b_{l_{n-2}}\cdots).
 \end{align*}

When $2=k=m< n=3$, we have $a_{h_{m-1}}=a_{h_1}=a_{t_2}=a_1$, $b_{t_1}=b_3$, $b_{l_1}=b_{l_{n-2}}=b_2$, and 
$$
f_{1}=g_{2}=(a_{1}b_{3}a_{2}b_{1}\cdots),\  g_{1}=(b_{2}a_{1}b_{3}\cdots).
 $$

When $3=k=m<n$, we have $\{a_{h_{m-1}}, a_{h_1}\}=\{a_1,a_2\}$, and  
$$f_{2}=g_{ \lceil n/2 \rceil}=(a_{h_{2}}b_{t_1}a_{3}b_{t_{2}}\cdots),\ f_1=(a_{h_{1}}b_{t_1}a_{h_{m-1}}\cdots).$$
  
Consequently, $$\big(\{f_1,f_2,\ldots, f_{\lceil m/2 \rceil}\},\{g_1,g_2,\ldots, g_{\lceil n/2 \rceil}\}\big)$$ is an extendible pair in $\Pi$.
 
{\bf Case 3.}\quad $4\leq k=m=n.$
 Suppose the rotation at vertex $b_{m}$   induced by $\Pi$ is 
 $$\Pi_{b_m}=[a_{i_1}a_{i_2}\cdots a_{i_{m-1}}].$$
For any vertex $b_x\neq b_m$, the rotation at $b_x$ has the  form
 $$\Pi_{b_x}=[a_{j_1}a_ma_{j_2}\cdots].$$
 Thus, after relabeling the vertices if necessary, we may assume that the rotations at $b_m$ and some vertex $b_w$ satisfy
 \begin{equation}
\label{R1}
\Pi_{b_m}=[a_{1}a_{2}\cdots a_{{m-1}}],\quad \Pi_{b_{w}}=[a_1a_{m}a_{q}\cdots].
\end{equation}


\begin{figure}[h]

\setlength{\unitlength}{1.5mm}
  \begin{center}
\begin{tikzpicture}[]

\def\e {1.5};
\def\f{4};
  \coordinate [label={[label distance=0mm]60:{\small$a_1$}}](a1) at (\f,4);
   \coordinate [label={[label distance=0mm]-10:{\small$b_{m}$}}] (bm) at (\f,1.5);
   \coordinate [label={[label distance=0mm]110:{\small$b_{w}$}}]  (bw) at ({\f+3}, {4}); 
   \coordinate [label={[label distance=0mm]90:{\small $a_{q+2}$}}](a2p+1) at (\f+4, 1.5);

\coordinate [label={[label distance=0mm]-10:{\small$a_{Y}$}}]  (aY) at (\f-4, 1.5);

\coordinate [label={[label distance=0mm]0:{\small $b^{'}$}}](b') at (({\f+1.5*cos(-30)}, {4+1.5*sin(-30)});
\draw[thick] (a1)--(b');
\filldraw [draw=black, fill=red](b') circle (2pt);

\foreach \i in {1,2,3} {
\node at ({\f+0.8*cos(-80+\i*10)}, {4+0.8*sin(-80+\i*10)}){$\cdot$};
}

\coordinate [label={[label distance=0mm]0:{\small $a_{m}$}}](am) at ({\f+3}, {2.5});
\coordinate [label={[label distance=0mm]0:{\small $a_{q}$}}](aq) at ({\f+3+1.5*cos(-45)}, {4+1.5*sin(-45)});
\draw[thick] (a1)--(bw)--(am)(bw)--(aq);
\fill[blue!20] (bw)--({\f+3+1.4*cos(-45)}, {4+1.4*sin(-45)})--({\f+3}, {2.6})--(bw);
\node at ({\f+3+1*cos(-67.5)}, {4+1*sin(-67.5)}){$f_{1}$};


 \draw[thick] (a1)--(bm)--(aY)(bm)--(a2p+1);
\filldraw [draw=black, fill=red](aq) circle (3pt)(am) circle (3pt);

 \coordinate [label={[label distance=0mm]0:{\small $b_{Y}$}}](bY) at ({\f+4}, {0});
\draw[thick] (bY)--(a2p+1) ;
 
\def\a{60};
\def\b{30};
\foreach \i in {3,2} {
\coordinate [label={[label distance=0mm]\a+\b*\i:{\small $a_{{\i}}$}}](a\i) at ({\f+1.5*cos(\a+\b*\i)}, {1.5+1.5*sin(\a+\b*\i)});
\draw[thick] (bm)--(a\i) ;
\filldraw [draw=black, fill=red](a\i) circle (3pt);
}

\coordinate [label={[label distance=0mm]-90:{\small $a_{q+1}$}}](a2p) at ({\f+1.5*cos(-90)}, {1.5+1.5*sin(-90)});
\draw[thick] (bm)--(a2p) ;\draw[thick,dotted] (bY)--(a2p);
 \fill[red!20] (bm)--(a2p)--(bY)--(a2p+1)--(bm);
 \filldraw [draw=black, fill=red](a2p) circle (3pt);
  \filldraw [draw=black, fill=blue](bY) circle (3pt);

\coordinate [label={[label distance=0mm]-120:{\small $a_{q}$}}](a2p-1) at ({\f+1.5*cos(-120)}, {1.5+1.5*sin(-120)});
\draw[thick] (bm)--(a2p-1) ;
\filldraw [draw=black, fill=red](a2p-1) circle (3pt);
\coordinate [label={[label distance=0mm]-150:{\small $a_{q-1}$}}](a2p-2) at ({\f+1.5*cos(-150)}, {1.5+1.5*sin(-150)});
\draw[thick] (bm)--(a2p-2) ;
\filldraw [draw=black, fill=red](a2p-2) circle (3pt);

\node at({\f+1*cos(70)}, {1.5+1*sin(70)}){$\cdot$};
\node at ({\f+1*cos(80)}, {1.5+1*sin(80)}){$\cdot$};
\node at ({\f+1*cos(75)}, {1.5+1*sin(75)}){$\cdot$};

\coordinate [label={[label distance=0mm]0:{\small $a_{q+3}$}}](a2i) at ({\f+1.5*cos(30)}, {1.5+1.5*sin(30)});
\coordinate [label={[label distance=0mm]0:{\small $a_{q+4}$}}](a2i+1) at ({\f+1.5*cos(60)}, {1.5+1.5*sin(60)});
\draw[thick] (a2i+1)--(bm)--(a2i);
\filldraw [draw=black, fill=red] (a2i) circle (3pt)(a2i+1) circle (3pt);
\fill[blue!20] (bm)--({\f+1.4*cos(30)}, {1.5+1.4*sin(30)})--({\f+1.4*cos(60)}, {1.5+1.4*sin(60)})--(bm);
\node at ({\f+1*cos(45)}, {1.5+1*sin(45)}){$f_{2}$};

\node at({\f+1*cos(160)}, {1.5+1*sin(160)}){$\cdot$};
\node at ({\f+1*cos(165)}, {1.5+1*sin(165)}){$\cdot$};
\node at ({\f+1*cos(170)}, {1.5+1*sin(170)}){$\cdot$};

\node at({\f+1*cos(190)}, {1.5+1*sin(190)}){$\cdot$};
\node at ({\f+1*cos(195)}, {1.5+1*sin(195)}){$\cdot$};
\node at ({\f+1*cos(200)}, {1.5+1*sin(200)}){$\cdot$};

\def\a{30};
\def\b{30};
\foreach \i in {1,2,3,4} {
\coordinate [label={[label distance=0mm]\a+\b*\i:{\small $b_{k_{\i}}$}}](bk\i) at ({\f-4+1.5*cos(\a+\b*\i)}, {1.5+1.5*sin(\a+\b*\i)});
\draw[thick] (aY) --(bk\i) ;
\filldraw [draw=black, fill=blue](bk\i) circle (3pt);
}
\fill[red!20] (aY)--({\f-4+1.4*cos(60)}, {1.5+1.4*sin(60)})--({\f-4+1.4*cos(90)}, {1.5+1.4*sin(90)})--(aY);
\node at ({\f-4+1*cos(75)}, {1.5+1*sin(75)}){$g_{1}$};
\fill[red!20] (aY)--({\f-4+1.4*cos(120)}, {1.5+1.4*sin(120)})--({\f-4+1.4*cos(150)}, {1.5+1.4*sin(150)})--(aY);
\node at ({\f-4+1*cos(135)}, {1.5+1*sin(135)}){$g_{2}$};

\coordinate [label={[label distance=0mm]200:{\small $b_{k_{2t-1}}$}}](bk5) at ({\f-4+1.5*cos(200)}, {1.5+1.5*sin(200)});
\draw[thick] (aY) --(bk5) ;
\coordinate [label={[label distance=0mm]230:{\small $b_{k_{2t}}$}}](bk6) at ({\f-4+1.5*cos(230)}, {1.5+1.5*sin(230)});
\draw[thick] (aY) --(bk6);
\coordinate [label={[label distance=0mm]0:{\small $b_{k_{m-2}}$}}](bk7) at ({\f-4+1.5*cos(280)}, {1.5+1.5*sin(280)});
\draw[thick] (aY) --(bk7) ;
\fill[red!20] (aY)--({\f-4+1.4*cos(200)}, {1.5+1.4*sin(200)})--({\f-4+1.4*cos(230)}, {1.5+1.4*sin(230)})--(aY);
\node at ({\f-4+1*cos(215)}, {1.5+1*sin(215)}){$g_{t}$};

\filldraw [draw=black, fill=blue](bk5) circle (3pt)(bk6) circle (3pt)(bk7) circle (3pt);
\foreach \i in {1,2,3} {
\node at ({\f-4+1*cos(160+\i*10)}, {1.5+1*sin(160+\i*10)}){$\cdot$};
}
\filldraw [draw=black, fill=blue](bk5) circle (3pt)(bk6) circle (3pt)(bk7) circle (3pt);
\foreach \i in {1,2,3} {
\node at ({\f-4+1*cos(235+\i*10)}, {1.5+1*sin(235+\i*10)}){$\cdot$};
}

\node at ({\f+2.1}, {0.6}){{\small $f_{\lceil m/2 \rceil}=g_{\lceil n/2 \rceil}$}};

 \coordinate [label={[label distance=0mm]-10:{\small$b_{m}$}}] (bm) at (\f,1.5);
 
\filldraw [draw=black, fill=blue] (bm) circle (3pt)(bw) circle (3pt);
\filldraw [draw=black, fill=red](a1)circle (3pt) (aY) circle (3pt)(a2p+1) circle (3pt);
\end{tikzpicture}
\end{center}

\caption{{\small The induced rotations at $b_m$, $a_{1}$ and $a_{Y}$ for $4\leq k=m=n$}}\label{F3}
 \end{figure}
%
%

Clearly,  there exists a face $(a_{q+1}b_m a_{q+2}b_{Y}\cdots)$ in $\Pi$, where the index of $a$ exceeds $m-1$, it is taken by modulo $m-1$. Note that $a_Y$ may be equal  $a_1$. Furthermore, 
 suppose the rotation at $a_Y$ is 
 $$\Pi_{a_{Y}}=[b_{k_1}b_{k_2}\cdots b_{k_{m-2}}b_m].$$
The rotations at vertices $b_m, a_1,b_w$ and $a_Y$ in this case are shown in Figure \ref{F3}. 
 
 Now,   We choose the following faces  (when the index of $a$ exceeds $m-1$, it is taken by modulo $m-1$): 
\begin{align*}
 f_{\lceil m/2 \rceil}=g_{\lceil m/2 \rceil}&=(a_{q+1}b_m a_{q+2}b_{Y}\cdots),\\
  f_{1}&= (a_mb_{w} a_q\cdots), \\
 f_{i}&=(a_{q+2i-1}b_{m}a_{q+2i}\cdots)  ~\text{for}~2\leq i \lfloor (m-2)/2\rfloor,\\
  g_{t}&=(b_{k_{2t-1}}a_{Y}b_{k_{2t}}\cdots)  ~\text{for}~1\leq t \leq \lfloor (m-2)/2\rfloor, \\
 \intertext{when $m=n$ is odd, }
 f_{(m-1)/2}&= (a_{q-2}b_{m} a_{q-1}\cdots),\\
 g_{(m-1)/2}&= (b_{k_{m-3}}a_{Y} b_{k_{m-2}}\cdots).\\
 \end{align*}
Consequently,  we obtain an extendible pair $\{\mathcal{F}_A, \mathcal{F}_B\}$, where
$$\mathcal{F}_A=\{f_1,f_2,\ldots,f_{\lceil m/2 \rceil}\}, ~~\mathcal{F}_B=\{g_1,g_2,\ldots,g_{\lceil m/2 \rceil}\}.$$
 \proofend

\section{Extend embedding}\label{sec3}
 
A \emph{crosscap} is depicted as a hole marked with a $\sim$ inside,  where antipodal points on its boundary are identified. 
 A \emph{handle} is drawn as a pair of circular holes (called its \emph{ends}) ,
 whose boundaries are identified with each other. 
If the identification reverses orientation, the handle is called a \emph {prohandle}; if the orientation is preserved, it is called an \emph{antihandle}.

Consider an embedding $\Pi$ of a graph $G$ and two faces $f_1, f_2\in \mathcal{F} (\Pi)$. Suppose

$$
f_1=( x_1\cdots x_2\cdots x_t \cdots), f_2=(y_1\cdots y_2\cdots y_t \cdots),
$$ 
where vertices may be repeated in  the sequence $(x_1,x_2,\ldots,x_t)$ and $(y_1,y_2,\ldots,y_t)$. 
Assume that $x_iy_i\notin E(G)$ for each $i$ with $1\leq i\leq t$. 

\begin{defn} 
\emph{Inserting a handle $H$ within $(f_1,f_2)$ carrying the edge set
$\{x_iy_i \mid 1\leq i\leq t\}$} is performed as follows (see  Figure \ref{F6}):
\begin{enumerate}
    \item Remove a small open disk,  called an \emph{end},  from the interior of each face $f_1$ and $f_2$. Denote  the boundaries of these disks by $C_l$ and $C_r$, respectively, and assign an orientation to each boundary. 
  \item Glue $C_l$ and $C_r$ according to the assigned orientations so that the $t$ edges $\{x_iy_i\mid 1\leq i\leq t\}$ pass through the handle $H$ without  crossing each other. 
\end{enumerate}
 \end{defn}
 
\begin{figure}[h]
\begin{center}
\setlength{\unitlength}{1.5mm}
  \begin{tikzpicture}[scale=0.8,bezier bounding box,
    whitenode/.style={circle, draw=black, fill=red, thick,minimum size=0.5mm, inner sep =2.5pt},  
    smallnode/.style={circle, draw=black, fill=black!20, thick,minimum size=0.5mm, inner sep =2pt}, 
     bignode/.style={circle, draw=black, fill=black!20, thick,minimum size=0.5mm, inner sep =3pt},
    bluenode/.style={circle, draw=black, fill=blue, thick,minimum size=0.5mm, inner sep =2.5pt},
    springnode/.style={circle, draw=black, fill=black, minimum size=0.5mm, inner sep =1pt},  
    black_thick/.style={line width=1.3pt},  
    blackedge/.style={line width=0.5pt},]

    \coordinate (H1) at (6,11);
 \coordinate (H2) at (12,11);
  \draw[thick] (H1) circle (2cm);
 \draw[thick,blue!40] (H1) circle (0.8cm);
  \draw[thick] (H2) circle (2cm);
 \draw[thick,blue!40] (H2) circle (0.8cm);
  \node at ($(H1)+(0,1.5)$) {$f_1$};
  \coordinate [label=left:$x_1$]  (w1) at ($(H1)+(180:2)$);
   \coordinate [label=above:$x_{2}$]  (w2) at ($(H1)+(120:2)$);
  \coordinate [label=above:$x_{3}$]  (w3) at ($(H1)+(60:2)$);
 \coordinate [label=below:$x_{t}$]  (w4) at ($(H1)+(-120:2)$);
  \coordinate [label=below:$x_{{t-1}}$]  (w5) at ($(H1)+(-60:2)$);
  
  \coordinate []  (e1) at ($(H1)+(180:0.8)$);
   \coordinate []  (e2) at ($(H1)+(120:0.8)$);
  \coordinate []  (e3) at ($(H1)+(60:0.8)$);
 \coordinate []  (e4) at ($(H1)+(-120:0.8)$);
  \coordinate []  (e5) at ($(H1)+(-60:0.8)$);

  \node at ($(H1)+(1.4,1)$) {$\vdots$}; \node at ($(H1)+(1.4,-0.5)$) {$\vdots$};
\draw [style=bluenode]  (w1) circle (4pt)(w2)circle(4pt)
(w3)circle(4pt)(w4)circle(4pt) (w5)circle(4pt);
    \draw(w1)--(e1)(w2)--(e2)(w3)--(e3)(w4)--(e4)(w5)--(e5);

 \node at ($(H2)+(0,1.5)$) {$f_2$};
  \coordinate [label=left:$y_1$]  (v1) at  ($(H2)+(180:2)$);
   \coordinate [label=above:$y_{i_2}$]  (v2) at  ($(H2)+(120:2)$);
  \coordinate [label=above:$y_{i_3}$]  (v3) at  ($(H2)+(60:2)$);
 \coordinate [label=below:$y_{i_{t}}$]  (v4) at  ($(H2)+(-120:2)$);
  \coordinate [label=below:$y_{i_{t-1}}$]  (v5) at  ($(H2)+(-60:2)$);
  
  \coordinate []  (e1) at ($(H2)+(180:0.8)$);
   \coordinate []  (e2) at ($(H2)+(120:0.8)$);
  \coordinate []  (e3) at ($(H2)+(60:0.8)$);
 \coordinate []  (e4) at ($(H2)+(-120:0.8)$);
  \coordinate []  (e5) at ($(H2)+(-60:0.8)$);
 
 \node at  ($(H2)+(1.4,1)$) {$\vdots$}; \node at ($(H2)+(1.4,-0.5)$)  {$\vdots$};
\draw(v1)--(e1);  \draw(v2)--(e2); \draw(v3)--(e3); \draw(v4)--(e4);
\draw(v5)--(e5); 
\draw [style=whitenode]  (v1) circle (4pt)(v2)circle(4pt)
(v3)circle(4pt)(v4)circle(4pt) (v5)circle(4pt) ;
\fill[blue!10] (H1) circle (0.8cm);  \fill[blue!10] (H2) circle (0.8cm);
 \node at (H2){$C_r$};  \node at (H1){$C_l$};
 \end{tikzpicture}
\end{center}

	\caption
	{Inserting a handle $H$ within $(f_1,f_2)$ carrying edge set $\{x_iy_i \mid 1\leq i\leq t\}$. $H$ is an antihandle if the ordered indices $\langle i_2,i_3, \cdots,i_t\rangle=\langle t,t-1, t-2,\cdots, 2\rangle$;  it is a prohandle if $\langle i_2,i_3, \cdots,i_t\rangle=\langle 2,3,\cdots, t\rangle$.}
	\label{F6}
\end{figure}
 
 In particular, if $f_1=(b_1\cdots b_2\cdots)$ and $f_2=(a_1\cdots a_2\cdots)$, the handle inserted within $(f_1,f_2)$ that carries the edge set $E(\{a_1,a_2\},\{b_1,b_2\})$ is referred to as an \emph{$X$-type handle}. It is straightforward to verify that an $X$-type handle may  either a prohandle or an antihandle.
 
\begin{figure}[h]
\setlength{\unitlength}{1.5mm}
  \begin{center}
\begin{tikzpicture}[]
    \coordinate (H1) at (2,1.5);
   \draw[very thick] (H1) circle (1.8cm);
 \draw[very thick] (H1) circle (0.8cm);
 \fill[blue!20] (H1) circle (0.8cm);
  \coordinate (H2) at (7.5,1.5);
  \coordinate (H3) at (9.6,1.5);
  \draw[thick] (H2) circle (0.8cm);
  \fill[blue!20] (H2) circle (0.8cm);
   \draw[very thick] (H3) circle (0.8cm);
  \fill[red!10] (H3) circle (0.8cm);

     \coordinate [label=left:{\small $b_1$}] (b1) at ({2+1.8*cos(180)},{1.5+1.8*sin(180)});
      \coordinate [label=right:{\small $b_2$}] (b2) at ({2+1.8*cos(0)},{1.5+1.8*sin(0)});
       \coordinate [label={[label distance=-0.40mm]0:{\small $1$}}] (c1) at ({2+0.8*cos(130)},{1.5+0.8*sin(130)});
        \coordinate  [label={[label distance=-0.40mm]0:{\small $2$}}] (c2) at ({2+0.8*cos(180)},{1.5+0.8*sin(180)});
        \coordinate [label={[label distance=-0.40mm]0:{\small $3$}}] (c3) at ({2+0.8*cos(230)},{1.5+0.8*sin(230)});
       \coordinate [label={[label distance=-0.40mm]180:{\small $6$}}] (c4) at ({2+0.8*cos(50)},{1.5+0.8*sin(50)});
        \coordinate [label={[label distance=-0.40mm]180:{\small $5$}}] (c5) at ({2+0.8*cos(0)},{1.5+0.8*sin(0)});
        \coordinate [label={[label distance=-0.40mm]177:{\small $4$}}] (c6) at ({2+0.8*cos(-50)},{1.5+0.8*sin(-50)});
      
      \coordinate [label={[label distance=-0.40mm]-90:{\small $6$}}] (d1) at ({7.5+0.8*cos(90)},{1.5+0.8*sin(90)});
        \coordinate [label={[label distance=-0.50mm]0:{\small $1$}}] (d2) at ({7.5+0.8*cos(180)},{1.5+0.8*sin(180)});
       \coordinate [label={[label distance=-0.40mm]90:{\small $2$}}] (d3) at ({7.5+0.8*cos(260)},{1.5+0.8*sin(260)});
        \coordinate [label={[label distance=-0.40mm]180:{\small $3$}}] (d4) at ({7.5+0.8*cos(320)},{1.5+0.8*sin(320)});
       \coordinate [label={[label distance=-0.40mm]180:{\small $4$}}] (d5) at ({7.5+0.8*cos(0)},{1.5+0.8*sin(0)});
       \coordinate [label={[label distance=-0.40mm]180:{\small $5$}}] (d6) at ({7.5+0.8*cos(30)},{1.5+0.8*sin(30)});
       \def \a{9.6};
        \coordinate [label={[label distance=-0.40mm]2:{\small $\alpha$}}] (e4) at ({\a+0.8*cos(210)},{1.5+0.8*sin(210)});
       \coordinate [label={[label distance=-0.40mm]0:{\small $\beta$}}] (e5) at ({\a+0.8*cos(180)},{1.5+0.8*sin(180)});
       \coordinate [label={[label distance=-0.40mm]0:{\small $\gamma$}}] (e6) at ({\a+0.8*cos(150)},{1.5+0.8*sin(150)});
        \coordinate [label={[label distance=-0.40mm]270:{\small $\delta$}}] (e1) at ({\a+0.8*cos(90)},{1.5+0.8*sin(90)});
              \coordinate [label={[label distance=-0.40mm]88:{\small $\delta$}}](e2) at ({\a+0.8*cos(260)},{1.5+0.8*sin(260)});
       \coordinate [label={[label distance=-0.40mm]90:{\small $\gamma$}}] (e3) at ({\a+0.8*cos(300)},{1.5+0.8*sin(300)});
       \coordinate [label={[label distance=-0.40mm]-90:{\small $\alpha$}}] (e7) at ({\a+0.8*cos(50)},{1.5+0.8*sin(50)});
       \coordinate [label={[label distance=-0.40mm]180:{\small $\beta$}}](e8) at ({\a+0.8*cos(-10)},{1.5+0.8*sin(-10)});
      \coordinate [label=above:{\small $b_3$}] (b3) at (9.5,3);
           \coordinate [label=left:{\small $a_1$}] (a1) at (6, 3);
             \coordinate [label=below:{\small $a_2$}] (a2) at (9.5, 0);
               \coordinate [label=right:{\small $a_3$}] (a3) at (12, 3);
                \coordinate [label=left:] (b4) at (6, 0);
                 \coordinate [label=left:] (b5) at (12, 0);
           
\draw [thick](d2)--(a1) --(d1);
\draw [thick](d3)--(a2);
\draw [thick,dashed](d4)--(e4);
\draw [thick,dashed](d5)--(e5);
\draw [thick,dashed](d6)--(e6);
\draw [thick,dashed](b3)--(e1);
\draw [thick,dashed](a2)--(e2);
\draw [thick,dashed](a2)--(e3);
\draw [thick,dashed](a3)--(e7);
\draw [thick,dashed](a3)--(e8);
\draw [thick](b3) --(a1)(b3) --(a3);\draw [thick,dotted](a3)--(b5)--(a2);\draw [thick,dotted](a1)--(b4)--(a2);
    \foreach \i in {1,2,3} {
\filldraw [draw=black, fill=blue](b\i) circle (3pt);
\filldraw [draw=black, fill=red](a\i) circle (3pt);
}
 \foreach \i in {1,2,3} {
\draw [thick](b1) --(c\i);}
 \foreach \i in {4,5,6} {
\draw [thick](b2) --(c\i);
}
 \fill[blue!10] (H1) circle (0.8cm);  \fill[blue!10] (H2) circle (0.8cm);
  \node at (H2){$C_r$};  \node at (H1){$C_l$};
 \node at (H3){$\sim$};  
\node at (2,3){$f_1$}; 
\node at (8.5,2.5){$f_2$};  
\node at (11,1){$f_3$};  

 \coordinate [label={[label distance=0mm]-90:{\small $6$}}] (d1) at ({7.5+0.8*cos(90)},{1.5+0.8*sin(90)});
        \coordinate [label={[label distance=0mm]0:{\small $1$}}] (d2) at ({7.5+0.8*cos(180)},{1.5+0.8*sin(180)});
       \coordinate [label={[label distance=0mm]90:{\small $2$}}] (d3) at ({7.5+0.8*cos(260)},{1.5+0.8*sin(260)});
        \coordinate [label={[label distance=0mm]180:{\small $3$}}] (d4) at ({7.5+0.8*cos(320)},{1.5+0.8*sin(320)});
       \coordinate [label={[label distance=0mm]180:{\small $4$}}] (d5) at ({7.5+0.8*cos(0)},{1.5+0.8*sin(0)});
       \coordinate [label={[label distance=0mm]180:{\small $5$}}] (d6) at ({7.5+0.8*cos(30)},{1.5+0.8*sin(30)});
    \coordinate [label={[label distance=0mm]0:{\small $1$}}] (c1) at ({2+0.8*cos(130)},{1.5+0.8*sin(130)});
        \coordinate  [label={[label distance=0mm]0:{\small $2$}}] (c2) at ({2+0.8*cos(180)},{1.5+0.8*sin(180)});
        \coordinate [label={[label distance=0mm]0:{\small $3$}}] (c3) at ({2+0.8*cos(230)},{1.5+0.8*sin(230)});
       \coordinate [label={[label distance=0mm]180:{\small $6$}}] (c4) at ({2+0.8*cos(50)},{1.5+0.8*sin(50)});
        \coordinate [label={[label distance=0mm]180:{\small $5$}}] (c5) at ({2+0.8*cos(0)},{1.5+0.8*sin(0)});
        \coordinate [label={[label distance=0mm]177:{\small $4$}}] (c6) at ({2+0.8*cos(-50)},{1.5+0.8*sin(-50)});    
     \foreach \i in {4,5} {
\filldraw [draw=black, fill=blue](b\i) circle (2pt);
}
 \end{tikzpicture}
\end{center}
\caption{\small Inserting a handle plus a crosscap within $\langle f_1,f_2, f_3\rangle$ carrying the edge set
$E(\{b_1, b_2\}, \{a_1,a_2,a_3\})$}.
	\label{F7}
\end{figure}

  \begin{defn} 
Let $f_1=(b_1\cdots b_2\cdots)$, $f_2=(a_1\cdots b_3 a_2\cdots)$ and $f_3=(b_3a_2\cdots a_3\cdots)$ be three faces in $\Pi$.
\emph{Inserting a handle $H$ plus a crosscap within $\langle f_1,f_2, f_3\rangle$ carrying the edge set
$E(\{b_1, b_2\}, \{a_1,a_2,a_3\})$} is performed by firstly inserting the handle $H$ within $(f_1,f_2)$ carrying the edge set $E(\{b_1, b_2\}, \{a_1,a_2,a_3\})$, and then deleting the edge $b_2a_2$ and adding a crosscap  in the vicinity of $b_{3}$  so that the edge set
 $\{a_{2}b_{2},a_2b_3, a_{3}b_{2}, a_{3}b_{1}\}$ pass through this crosscap.  
See Figure \ref{F7}. 
\end{defn}
 
\begin{defn} 
Let $f_1=(\cdots b_1 a_2 b_2\cdots)$ and $f_2=(b_xa_2b_y\cdots a_1\cdots)$ be two distinct faces in $\Pi$, where the rotation at $a_2$ induced by $\Pi$ is 
$$\Pi_{a_2}=[\cdots b_1b_2b_3\cdots b_k b_xb_y\cdots].$$
The operation of \emph{inserting a crosscap $C$ near $a_2$ between $\langle f_1,f_2\rangle$ carrying the edge set $\{b_1a_1, b_2a_1\}$} is performed by firstly removing the edges $a_2b_2, a_2b_3, \cdots, a_2b_k, a_2b_x$, and then adding a crosscap in the vicinity of $a_2$ so that the edge set $\{a_2b_2, a_2b_3, \cdots, a_2b_k, a_2b_x\}\cup \{b_1a_1, b_2a_1\}$ pass through this crosscap. See Figure \ref{F8}. 
\end{defn}

\begin{figure}[h]
\setlength{\unitlength}{1.5mm}
  \begin{center}
\begin{tikzpicture}[]

   \def\r{1.5};
     \coordinate [label={[label distance=0mm]-90:{\small $a_2$}}] (a2) at (-3,1);
      \coordinate [label={[label distance=0mm]180:{\small $b_{1}$}}] (b1) at (-3-\r,1);

     \coordinate [label={[label distance=0mm]90:{\small $b_{2}$}}] (b2) at ({-3+\r*cos(130)},{1+\r*sin(130)});
\coordinate [label={[label distance=0mm]90:{\small $b_{3}$}}] (b3) at ({-3+\r*cos(110)},{1+\r*sin(110)});
\coordinate [label={[label distance=0mm]90:{\small $b_{k}$}}] (bk) at ({-3+\r*cos(70)},{1+\r*sin(70)});
\coordinate [label={[label distance=0mm]90:{\small $b_{x}$}}] (bx) at ({-3+\r*cos(50)},{1+\r*sin(50)});
\coordinate [label={[label distance=0mm]0:{\small $b_{y}$}}] (by) at ({-3+\r*cos(0)},{1+\r*sin(0)});
      \coordinate [label={[label distance=0mm]90:{\small $$}}] (a1) at ({-3-\r+\r*cos(130)},{1+\r*sin(130)});
       \coordinate [label={[label distance=0mm]90:{\small $a_{1}$}}] (a3) at ({-3+\r+\r*cos(50)},{1+\r*sin(50)});
        \fill[gray!10] (a1) --(b1)--(a2)--(b2)--(a1); \fill[gray!10] (a3) --(by)--(a2)--(bx)--(a3);
       
  \foreach \i in {1,2,3} {
\filldraw [draw=black, fill=red](a\i) circle (3pt);
}

 \foreach \i in {1,2,3,k,x,y} {
\filldraw [draw=black, fill=blue](b\i) circle (3pt);
}

\draw[thick](a2)--(b1)(a2)--(b2)(a2)--(b3)(a2)--(bk)(a2)--(bx)(a2)--(by);
  \draw[dotted, thick](a1)--(b1)(a1)--(b2)(a3)--(bx)(a3)--(by);
   
    \node at (1, 3.5){Inserting a crosscap near $a_2$};

    \node at (0.7, 1.5){$\Longrightarrow$};
 \node at (-4.3, 1.6){$f_1$};
     \node at (-2, 1.6){$f_2$};
    
    \foreach \i in {1,2,3} {
\node at ({-3+0.55*cos(255+8*\i)},{2.5+0.55*sin(255+8*\i)}){$\cdot$}; 
}

 \foreach \i in {1,2,3} {
\node at ({-3+0.55*cos(250+10*\i)},{1+0.55*sin(250+10*\i)}){$\cdot$}; 
}

   \coordinate (H2) at (5,2);
   \fill[red!10] (H2) circle (0.4cm);
  
 \draw[thin] (H2) circle (0.4cm);
  
  \node at (H2) {$\sim$};
 
  \def\r{2};
     \coordinate [label={[label distance=0mm]-90:{\small $a_2$}}] (a2) at (5,1);
      \coordinate [label={[label distance=0mm]180:{\small $b_{1}$}}] (b1) at (5-\r,1);

     \coordinate [label={[label distance=0mm]90:{\small $b_{2}$}}] (b2) at ({5+\r*cos(130)},{1+\r*sin(130)});
\coordinate [label={[label distance=0mm]90:{\small $b_{3}$}}] (b3) at ({5+\r*cos(110)},{1+\r*sin(110)});
\coordinate [label={[label distance=0mm]90:{\small $b_{k}$}}] (bk) at ({5+\r*cos(70)},{1+\r*sin(70)});
\coordinate [label={[label distance=0mm]90:{\small $b_{x}$}}] (bx) at ({5+\r*cos(50)},{1+\r*sin(50)});
\coordinate [label={[label distance=0mm]0:{\small $b_{y}$}}] (by) at ({5+\r*cos(0)},{1+\r*sin(0)});
      \coordinate [label={[label distance=0mm]90:{\small $$}}] (a1) at ({5-\r+\r*cos(130)},{1+\r*sin(130)});
       \coordinate [label={[label distance=0mm]90:{\small $a_{1}$}}] (a3) at ({5+\r+\r*cos(50)},{1+\r*sin(50)});
        \fill[gray!10] (a1) --(b1)--(a2)--(b2)--(a1); \fill[gray!10] (a3) --(by)--(a2)--(bx)--(a3);

\draw [thick](a2)--(b1)(a2)--(by);
\draw [thick,dotted](a1)--(b1)(a1)--(b2)(a3)--(bx)(a3)--(by);
 \foreach \i in {1,2} {
\draw [draw=black,densely dashed](a2)--({5+0.4*cos(200+20*\i)},{2+0.4*sin(200+20*\i)});
}
\foreach \i in {1,2} {
\draw [draw=black,densely dashed](a2)--({5+0.4*cos(270+20*\i)},{2+0.4*sin(270+20*\i)});
}

\foreach \i in {1,2,3} {
\node at ({5+0.55*cos(255+6*\i)},{2+0.55*sin(255+6*\i)}){$\cdot$}; 
}
\draw [draw=black,thick](b1)--({5+0.4*cos(210)},{2+0.4*sin(210)});
\draw [draw=black,thick](b2)--({5+0.4*cos(185)},{2+0.4*sin(185)});
\draw [draw=black,densely dashed](b2)--({5+0.4*cos(160)},{2+0.4*sin(160)});
\draw [draw=black,densely dashed](b3)--({5+0.4*cos(130)},{2+0.4*sin(130)});
\draw [draw=black, densely dashed](bk)--({5+0.4*cos(100)},{2+0.4*sin(100)});
\draw [draw=black,densely dashed](bx)--({5+0.4*cos(50)},{2+0.4*sin(50)});
\draw [draw=black,thick](a3)--({5+0.4*cos(20)},{2+0.4*sin(20)});
\draw [draw=black,thick](a3)--({5+0.4*cos(-10)},{2+0.4*sin(-10)});
\foreach \i in {1,2,3} {
\node at ({5+0.7*cos(110-7*\i)},{2+0.7*sin(110-6*\i)}){$\cdot$}; 
}
\foreach \i in {1,2,3} {
\node at ({5+0.55*cos(250+10*\i)},{1+0.55*sin(250+10*\i)}){$\cdot$}; 
}
 \foreach \i in {1,2,3} {
\filldraw [draw=black, fill=red](a\i) circle (3pt);
}

 \foreach \i in {1,2,3,k,x,y} {
\filldraw [draw=black, fill=blue](b\i) circle (3pt);
}
 \end{tikzpicture}
\end{center}
\caption{\small Inserting a crosscap near $a_2$ between $\langle f_1, f_2\rangle$ carrying  the edge set 
 $\{b_{1}a_{3}, b_{2}a_{3}\}$.}
	\label{F8}
\end{figure}
  
In the following of this section, let  $(\mathcal{F}_A, \mathcal{F}_B)$ be an extendible pair in the  embedding $\Pi$ of $G_{A,B}$, 
where
 \begin{align*}
 f_{\lceil m/2\rceil}=g_{\lceil n/2\rceil}&=(a_{m-1}\cdots b_{n-1} \cdots a_{m}\cdots  b_{n} \cdots),\\
f_i&=(a_{2i-1}\cdots b'_i a_{2i}\cdots b''_i) ~\text{for}~1\leq i\leq \lfloor (m-2)/2\rfloor, \\
g_j&=(b_{2j-1}\cdots a'_j b_{2j}\cdots a''_j)~\text{for}~1\leq j\leq  \lfloor (n-2)/2\rfloor, \\
\intertext{if $m$ is odd, }
 f_{ (m-1)/2}&=(a_{m-3}\cdots b'_{(m-1)/2} a_{m-2}\cdots b''_{(m-1)/2}),\\
\intertext{if $n$ is odd, }
 g_{ (n-1)/2}&=(b_{n-3}\cdots a'_{(n-1)/2}  b_{n-2}\cdots a''_{(n-1)/2}\cdots ).\\
\end{align*}
In particular, if $m=3$, then 
$$\mathcal{F}_A=\{f_1=(a_{1}\cdots b'_1a_{2}\cdots b''_1), \quad f_2=(a_{2}\cdots b_{n-1} \cdots a_{3}\cdots  b_{n} \cdots)\},$$ 
if $n=3$, then 
$$\mathcal{F}_B=\{g_1=(b_{1}\cdots a'_1b_{2}\cdots a''_1), \quad g_2=(a_{m-1}\cdots b_{2} \cdots a_{m}\cdots  b_{3} \cdots)\}.$$

In this and the next section, for any integer $x\geq 0$ we define 
 \begin{equation*}
	\tau_{x}=
	\begin{cases}
		\displaystyle 1, &~ \text{if}~ x ~ \text{is odd},\\[0.3cm]
		\displaystyle 0, &~ \text{if}~ x ~ \text{is even}.
	\end{cases}
\end{equation*}

\begin{defn}\label{defn3-3}
\emph{Inserting $G(2,2,2)$ in an embedding $\Pi$ of $G_{A,B}$}, the obtain embedding is denoted by $\Pi \oplus G(2,2,2)$, is constructed in the following four steps:
\begin{enumerate}
  \item[{\bf Step1.}] Adding the vertices $a_{m+1}$, $a_{m+2}$, $b_{n+1}$, $b_{n+2}$ and the edges $$E(\{a_{m-1}, a_{m}\}, \{b_{n+1},b_{n+2}\}) \cup E(\{b_{n-1}, b_n\}, \{a_{m+1},a_{m+2}\})$$ in $f_{\lceil m/2\rceil}$$=g_{\lceil n/2\rceil}$, see Figure \ref{F9}(1).
The vertices $a_{m+1}$, $a_{m+2}$, $b_{n+1}$, $b_{n+2}$ and edges 
$a_{m+1}b_{n+2}$, $a_{m+2}b_{n+1}$  are on the boundary of the crosscap.
This step introduces two new faces that will be used in the subsequent steps:
 $$g_{\lceil n/2\rceil+1}=(a_{m}b_{n+1}a_{m-1}b_{n+2}),~~f_{\lceil m/2\rceil+1}=(b_{n}a_{m+1}b_{n-1}a_{m+2}).$$
   
 \item[{\bf Step2.}] For each $i$ with $1\leq i\leq \lfloor (m-2)/2\rfloor-\tau_m$, inserting an $X$-type handle $H_i$ within $(f_{i}, g_{\lceil n/2\rceil+1})$ carrying the edge set $E(\{b_{n+1}, b_{n+2}\}, \{a_{2i-1}, a_{2i}\}).$ 
 For each $j$ with $1\leq j\leq \lfloor (n-2)/2\rfloor-\tau_n$, inserting an $X$-type handle $H'_j$ within $(g_{j}, f_{\lceil m/2\rceil+1})$ carrying the edge set $E(\{a_{m+1}, a_{m+2}\}, \{b_{2j-1},b_{2j}\}\}).$

  \item[{\bf Step3.}] If $m>3$ is odd, inserting a handle $H_{\lfloor (m-2)/2\rfloor}$ plus a crosscap $C$ within $\langle g_{\lceil n/2\rceil+1}$, $f_{\lfloor (m-2)/2\rfloor}$, $f_{(m-1)/2}\rangle$ carrying the edge set 
  $E(\{b_{n+1}, b_{n+2}\}$, $\{a_{m-4}, a_{m-3},a_{m-2}\}).$ 
 
  If $n>3$ is odd, inserting a handle $H'_{\lfloor (n-2)/2\rfloor}$ plus a crosscap $C'$ within $\langle f_{\lceil m/2\rceil+1}$, $g_{\lfloor (n-2)/2\rfloor}$, $g_{(n-1)/2}\rangle$ carrying the edge set 
  $E(\{a_{m+1}, a_{m+2}\}$, $\{b_{n-4}, b_{n-3},b_{n-2}\})$. 
  
 \item[{\bf Step4.}] If $m=3$, inserting a crosscap near vertex $a_2$ between $\langle g_{\lceil n/2\rceil+1}, f_1\rangle$ carrying 
  the edge set $\{b_{n+2}a_{1}, b_{n+1}a_{1}\}$.  
 If $n=3$,  inserting a crosscap near vertex $b_2$ between $\langle f_{\lceil m/2\rceil+1}, g_1\rangle$ carrying 
  the edge set $\{a_{m+2}b_{1}, a_{m+1}b_{1}\}$. See Figure \ref{F8}.
  
\end{enumerate}
\end{defn}
   
 \begin{figure}[h]
\setlength{\unitlength}{0.5mm}
\begin{center}
  \begin{tikzpicture}
 \begin{scope}[xshift=-4cm,yshift=0cm]
  \coordinate  [label={[label distance=-0.4mm]93:{\small $b_{n+1}$}}]  (bn+1) at ({9.5+0.8*cos(120)},{3.5+0.8*sin(120)});
  \coordinate [label={[label distance=0mm]180:{\small$b_{n+2}$}}]  (bn+2) at ({9.5+0.8*cos(165)},{3.5+0.8*sin(165)});
    \coordinate [label={[label distance=0mm]-80:{\small$b_{n+2}$}}]  (b''n+2) at ({9.5+0.8*cos(300)},{3.5+0.8*sin(300)});
  \coordinate [label={[label distance=0mm]0:{\small$b_{n+1}$}}]  (b''n+1) at ({9.5+0.8*cos(-15)},{3.5+0.8*sin(-15)});
   \coordinate [label={[label distance=-0.1mm]180:{\small$a_{m+1}$}}]  (am+1) at({9.5+0.8*cos(210)},{3.5+0.8*sin(210)});
   \coordinate [label={[label distance=0mm]-90:{\small$a_{m+2}$}}]  (am+2) at({9.5+0.8*cos(255)},{3.5+0.8*sin(255)});
   \coordinate [label={[label distance=0mm]0:{\small$a_{m+2}$}}]  (a''m+2) at ({9.5+0.8*cos(30)},{3.5+0.8*sin(30)});
   \coordinate [label={[label distance=-0.4mm]86:{\small$a_{m+1}$}}]  (a''m+1) at ({9.5+0.8*cos(75)},{3.5+0.8*sin(75)});
   \coordinate [label={[label distance=0mm]0:{\small$b_{n-1}$}}]  (bn-1) at (11.5,5);
  \coordinate [label={[label distance=0mm]180:{\small$b_{n}$}}]  (bn) at (7.5,2);
 \coordinate [label={[label distance=0mm]180:{\small$a_{m}$}}]  (am) at (7.5,5);
   \coordinate [label={[label distance=0mm]0:{\small$a_{m-1}$}}]  (am-1) at (11.5,2);  
   \draw (am)--(bn) 
   (am)--(bn-1)--(am-1)
   (am)--(bn+1)(am)--(bn+2)
    (am-1)--(b''n+1) (am-1)--(b''n+2) (bn-1)--(a''m+1)(bn-1)--(a''m+2) (bn)--(am+1)(bn)--(am+2) ;
      \draw[dashed, thick] (bn)--(am-1);
       \fill[fill=red!20] (9.5,3.5) circle (0.8cm) ;
       \node at (9.5,3.5) {$\sim$};
       \draw[thick]({9.5+0.8*cos(165)},{3.5+0.8*sin(165)}) arc (165:210:0.8);  \draw[densely dotted,  thick]({9.5+0.8*cos(210)},{3.5+0.8*sin(210)}) arc (210:255:0.8); 
 \draw[thick]({9.5+0.8*cos(255)},{3.5+0.8*sin(255)}) arc (255:300:0.8);  \draw[densely dotted, very thick]({9.5+0.8*cos(300)},{3.5+0.8*sin(300)})  arc (300:345:0.8); 
  \draw[thick]({9.5+0.8*cos(345)},{3.5+0.8*sin(345)}) arc (345:390:0.8);   
         \draw[densely dotted, very thick]({9.5+0.8*cos(30)},{3.5+0.8*sin(30)}) arc (30:75:0.8);   \draw[thick] ({9.5+0.8*cos(75)},{3.5+0.8*sin(75)}) arc (75:120:0.8);
            \draw[densely dotted, very thick]({9.5+0.8*cos(120)},{3.5+0.8*sin(120)}) arc (120:165:0.8); 
    \filldraw [draw=black,fill=blue] (bn) circle (3pt) (bn-1) circle (3pt) (bn+1) circle (3pt) (bn+2) circle (3pt) (b''n+1) circle (3pt) (b''n+2) circle (3pt);   
         \filldraw [draw=black,fill=red] (am) circle (3pt) (am-1) circle (3pt)  (a''m+1) circle (3pt) (a''m+2) circle (3pt)(am+1) circle (3pt) (am+2) circle (3pt);
      \draw [black] (am) circle (3pt) (am-1) circle (3pt)
      (am+1) circle (3pt) (am+2) circle (3pt) (a''m+1) circle (3pt) (a''m+2) circle (3pt);
\node at (9.5,1.5) {$(1)$};
\end{scope}
   \begin{scope}[xshift=10cm,yshift=0cm]
   \coordinate [label=above:{\small$a_{m+1}$}]  (am+1) at (2,4);
   \coordinate [label=right:{\small $a_{m+2}$}]  (am+2) at (3,3);
     \coordinate [label=right:{\small$b_{n-1}$}]  (bn-1) at (4,5);
  \coordinate [label=left:{\small$b_{n}$}]  (bn) at (1,2);
 \coordinate [label=above:{\small$a_{m}$}]  (am) at (1,5);
   \coordinate [label=right:{\small$a_{m-1}$}]  (am-1) at (4,2);  
         \draw (am)--(bn) 
   (am)--(bn-1)--(am-1)
 (bn)--(am+1)(bn)--(am+2) (am+2)--(bn-1)(am+1)--(bn-1);
  \fill[gray!20](am+1)--(bn)--(am+2)--(bn-1)--(am+1);
      \draw[dashed, thick] (bn)--(am-1);
    \filldraw [draw=black,fill=blue] (bn) circle (3pt) (bn-1) circle (3pt);   
         \filldraw [draw=black,fill=red] (am) circle (3pt) (am-1) circle (3pt)  (am+1) circle (3pt) (am+2) circle (3pt);
      \draw [black] (am) circle (3pt) (am-1) circle (3pt)
      (am+1) circle (3pt) (am+2) circle (3pt) ;
\node at (2.5,1.5) {$(2)$};\node at (2.8,3.5){{\small$f_{\lceil (m)/2\rceil+1}$}};
\end{scope}
\end{tikzpicture}
\caption{{\small (1)~adding two vertices $a_{m+1}, a_{m+2}$ in the face $f_{\lceil m/2\rceil}$;
 (2)~adding four vertices $a_{m+1}, a_{m+2}, b_{n+1}, b_{n+2}$ in the face $f_{\lceil m/2\rceil}$ through a crosscap.
}}\label{F9}
\end{center}
\end{figure}
    
 Similarly, we give the definition of $\Pi \oplus G(2,0,0)$ and  $\Pi \oplus G(0,2,0)$ by inserting $G(2,0,0)$ or $G(0,2,0)$ in  an embedding  $\Pi$ of $G_{A,B}$. By the symmetry, we only introduce the construction of $\Pi \oplus G(2,0,0)$.

 \begin{defn} \label{defn3-4}
 \emph{Inserting $G(2,0,0)$ in the embedding $\Pi$ of $G_{A, B}$}, denoted the obtained embedding by $\Pi \oplus G(2,0,0)$,  is constructed as follows. First,  add vertices $a_{m+1}$, $a_{m+2}$ into the face $f_{\lceil m/2\rceil}$ and connect all edges in $E(\{a_{m+1}, a_{m+2}\},\{b_n,b_{n-1}\})$, as shown in Figure \ref{F9}(2).
This step introduces a new face
  $$f_{\lceil m/2\rceil+1}=(b_{n}a_{m+1}b_{n-1}a_{m+2}).$$
 Second, add the edge set $E(\{a_{m+1}, a_{m+2}\},B\setminus \{b_n,b_{n-1}\})$ following  the procedure described in Step$2$ to Step$4$ of Definition \ref{defn3-3}.
 \end{defn}
  
   \begin{defn} \label{defn3-5}
 Let $\Pi$ be an extendible embedding of a bipartite graph. For $G(r,s,t)\in \{G(2,2,2), G(2,0,0), G(0,2,0)\}$, define $\Pi\oplus G(r,s,t)^t$ as the \emph{extend embedding} obtained from $\Pi$ by iteratively inserting $t$ copies of $G(r,s,t)$. That is 
 $$\Pi\oplus G(r,s,t)^t=\Pi\oplus G(r,s,t)\oplus \cdots \oplus G(r,s,t)$$
  with $t$ copies of $G(r,s,t)$.
   \end{defn}

 \begin{thm}\label{thm3}
Let $\Pi$ be  any embedding of $G(m,n,k)$ with $\mathcal{E}(\Pi)=\lceil f(m,n,k)\rceil$, $m,n\geq 2$,
 then 
\begin{align*}
\tilde{\gamma}(\Pi\oplus G(2,2,2))=&\lceil  f(m+2,n+2,k+2)\rceil,\\
\mathcal{E}(\Pi\oplus G(2,0,0))=&\lceil  f(m+2,n,k)\rceil, \\
 \mathcal{E}(\Pi\oplus G(0,2,0)) =&\lceil  f(m,n+2,k)\rceil.
\end{align*}
\end{thm}

\textbf{Proof.}
Since  when $(m,n,k)=(2,2,2), (3,3,3)$, $\lceil  f(m,n,k)\rceil =-1<0$, $G(m,n,k)$ do not have an embedding with Euler genus $\lceil f(m,n,k)\rceil$.  Hence, $G(m,n,k)\neq$ $(2,2,2)$, $(3,3,3)$. 

 By the definition \ref{defn3-3}, $\Pi\oplus G(2,2,2)$ is obtained from $\Pi$ by adding $1+\tau_m+\tau_n$ crosscap and $\lfloor (m-2)/2\rfloor+\lfloor (n-2)/2\rfloor$ handles, hence, we have 
\begin{eqnarray*}
 \tilde{\gamma}(\Pi\oplus G(2,2,2))&=&\mathcal{E}(\Pi)+1+2 \Big\lfloor \frac{m-2}{2}\Big\rfloor+2 \Big\lfloor \frac{n-2}{2} \Big\rfloor+\tau_m+\tau_n \\
 &=&\lceil  f(m,n,k)\rceil+m+n-3\\
 &=&\lceil f(m+2,n+2,k+2)\rceil.
\end{eqnarray*}

By the definition \ref{defn3-4}, $\Pi\oplus G(2,0,0)$ is obtained from $\Pi$ by adding $\tau_n$ crosscap and $\lfloor(n-2)/2 \rfloor$ handles, hence, we have 
$$
\mathcal{E}(\Pi\oplus G(2,0,0))=\mathcal{E}(\Pi)+2 \Big\lfloor \frac{n-2}{2}\Big\rfloor+\tau_n=\lceil  f(m,n,k)\rceil+n-2=\lceil f(m+2,n,k)\rceil.
$$
Similarly, we have 
$$\mathcal{E}(\Pi\oplus G(0,2,0))=\mathcal{E}(\Pi)+2 \Big\lfloor \frac{m-2}{2}\Big\rfloor+\tau_m=\lceil  f(m,n,k)\rceil+m-2=\lceil f(m,n+2,k)\rceil.$$
 \proofend

  \section{The proof of Theorem \ref{thm2}}\label{sec4}
 
 In this section, by the symmetry of the nearly complete bipartite graph, we assume that $k\leq m\leq n$. Since  $G(3,3,3)$ and $G(4,4,4)$ are planar graphs, we have 
$$\tilde{\gamma}(G(3,3,3))=\max\{  f(3,3,3), \,1\},\quad \tilde{\gamma}(G(4,4,4))=\max\{  f(4,4,4), \,1\}.$$
Thus, in what follows, we suppose that 
$$(m,n,k)\notin \{(3,3,3), (4,4,4), (5,4,4),(4,5,4), (5,5,5)\}.$$  
  
  \begin{figure}[h]
\begin{center}

\setlength{\unitlength}{1.5mm}
  \begin{tikzpicture}
  \def\a{-0.5};
  \def\b{6.5};
  \def\c{8};
  \def\d{14};
\coordinate [label=above:$b_1$]  (b11) at ({1},4);
\coordinate [label=below:$b_1$]  (b12) at ({1},0);
\coordinate [label=above:$b_2$]  (b21) at ({2},4);
\coordinate [label=below:$b_2$]  (b22) at ({2},0);
\coordinate [label=above:$b_3$]  (b31) at ({3},4);
\coordinate [label=below:$b_3$]  (b32) at ({3},0);
\coordinate [label=above:$b_4$]  (b41) at ({4},4);
\coordinate [label=below:$b_4$]  (b42) at ({4},0);
\coordinate [label=above:$b_5$]  (b51) at ({5},4);
\coordinate [label=below:$b_5$]  (b52) at ({5},0);

\coordinate [label=right:$a_2$]  (a2) at ({1},2);
\coordinate [label=right:$a_3$]  (a3) at ({2},2);
\coordinate [label=right:$a_4$]  (a4) at ({3},2);
\coordinate [label=right:$a_5$]  (a5) at ({4},2);
\coordinate [label=right:$a_1$]  (a1) at ({5},2);

\draw [thick]  (b11)--(a2) (b11)--(a3)(b12)--(a4) (b12)--(a5);
\draw [thick]  (b21)--(a3) (b21)--(a4)(b22)--(a5) (b22)--(a1);
\draw [thick]  (b31)--(a4) (b31)--(a5)(b32)--(a1) (b32)--(\b, 2)(\a, 2)--(a2);
\draw [thick]  (b41)--(a5) (b41)--(a1)(b42)--(\b, 1.6)(\a, 1.6)--(a2) (b42)--(\b, 1.2)(\a, 1.2)--(a3);
\draw [thick]  (b51)--(a1) (b51)--(\b, 3)(\a, 3)--(a2)(b52)--(\b, 0.8)(\a, 0.8)--(a3) (b52)--(\b, 0.4)(\a, 0.4)--(a4);

\draw [densely dashed, very thick,red]  (\a,0)--(\b,0)-- (\b,4)--(\a,4)--(\a,0);
\draw [densely dashed, very thick,red,->](\a, 2.2)--(\a, 2.4);
\draw [densely dashed, very thick,red,->](\b, 2.2)--(\b, 2.4);
\draw [densely dashed, very thick,red,->](2.4, 4)--(2.6, 4);
\draw [densely dashed, very thick,red,->](2.4, 0)--(2.6, 0);
 
\foreach \i in {1,2,3,4,5} {
\filldraw [draw=black,fill=red]  (a\i) circle (3pt);
}
\foreach \i in {11,21,31,41,51,12,22,32,42,52} {
\filldraw [draw=black,fill=blue]  (b\i) circle (3pt);
}

\filldraw [draw=black,fill=red!20] (10.5,1.8)circle (0.5cm);
\node at (10.5,1.8){$\sim$};
 \coordinate [label=above:$b_4$]  (b41) at ({\c+1},4);
\coordinate [label=right:$b_4$]  (b42) at ({\c+1},0);

\coordinate [label=below:$a_2$]  (a2) at ({\c+1.5},1);
\coordinate [label=above:$a_3$]  (a3) at ({\d-1.5},1);
\coordinate [label=below:$a_4$]  (a4) at ({\d-1.5},3);
\coordinate [label=above:$a_1$]  (a1) at ({\c+1.5},3);

\coordinate [label=below:$b_6$]  (b61) at (10.5,0);
\coordinate [label=above:$b_6$]  (b62) at (10.5,4);

\coordinate [label=above:$b_5$]  (b5) at (10.5,3);
\coordinate [label=left:$b_7$]  (b7) at ({\c+1.5},2);

\coordinate [label=above:$b_2$]  (b22) at ({\d-1},4);
\coordinate [label=below:$b_2$]  (b21) at ({\d-1},0);

\coordinate [label=left:$b_3$]  (b31) at (\c,3);
\coordinate [label=right:$b_3$]  (b32) at (\d,3);

\coordinate [label=left:$b_1$]  (b11) at (\c,1);
\coordinate [label=right:$b_1$]  (b12) at (\d,1);

\foreach \i in {1,4} {
\draw [thick]  (b5)--(a\i);
}
\foreach \i in {1,2} {
\draw [thick]  (b7)--(a\i);
}

\draw [densely dashed, very thick,red]  (\c,0)--(\d,0)-- (\d,4)--(\c,4)--(\c,0);
\draw [densely dashed, very thick,red,->](\c, 1.8)--(\c, 2);
\draw [densely dashed, very thick,red,->](\d, 1.8)--(\d, 2);
\draw [densely dashed, very thick,red,->](9.8, 4)--(10, 4);
\draw [densely dashed, very thick,red,->](9.8, 0)--(10, 0);

\draw [thick]  (b41)--(a1);
\draw [thick]  (b61)--(a3)(b61)--(a2)(b62)--(a1)(b62)--(a4);
\draw [thick]  (b21)--(a3)(a3)(b22)--(a4)(b22)--(\d,3.5)(\c,3.5)--(a1);

\draw [thick]  (b31)--(a1)(b31)--(a2)(b32)--(a4);
\draw [thick]  (b11)--(a2)(b12)--(a3)(b12)--(a4);
\draw [thick]  (b42)--(a2)(b42)--(\c, 0.6)(\d,0.6)--(a3);

\foreach \i in {1,2,3,4} {
\filldraw [draw=black,fill=red]  (a\i) circle (3pt);
}
\foreach \i in {11,12,21,22,31,32,41,42,5,7, 61,62} {
\filldraw [draw=black,fill=blue]  (b\i) circle (3pt);
}
\foreach \i in {1,2,3,4,5,6,7,8} {
\coordinate [label=left:]  (c\i) at ({10.5+0.5*cos(-22.5+45*\i)},{1.8+0.5*sin(-22.5+45*\i)});
}

\draw [densely dotted, thick]  (b5)--(c2)(b5)--(c3)(b7)--(c4)(b7)--(c5)(a2)--(c6)(a3)--(c7)(a3)--(c8)(a4)--(c1);

\node at (3,-1){$(1)$~$\Pi({5,5,5})$};\node at (11,-1){$(2)$~$\Pi({4,7,4})$};

\end{tikzpicture}

\caption{{\small $\Pi_{5,5,5}$ with $\mathcal{E}(\Pi(5,5,5)=f(5,5,5)=2$  and $\Pi_{4,7,4}$ with $\tilde{\gamma}(\Pi(4,7,4)=  f(4,7,4)=3$. The red rectangle represents an prohandle obtained by gluing the opposite sides of the rectangle along the directions indicated by the arrows.}
 }
\label{F10}
\end{center}
\end{figure}
  
To prove Theorem \ref{thm2}, it is necessary to identify a nearly complete bipartite graph $G(p,q,h)$ from Table \ref{table1} 
such that
  \begin{equation*}
G(m,n,k)= G(p, q, h) \oplus G(2,2,2)^{a} \oplus G(2,0,0)^{b}\oplus G(0,2,0)^{c},
\end{equation*}
where 
$$p+2a+2b=m,~ q+2a+2c=n,~ h+2a=k.$$

 The graph $G(p,q,h)$ is derived from $G(m,n,k)$ by first removing several copies of $G(2,0,0)$ or $G(0,2,0)$ to obtain $G(m', n', k)$,  where $k\leq m'\leq k+2$ and $m'\leq n'\leq m'+2$; then, by further  removing a number of copies of $G(2,2,2)$ to arrive at $G(p,q,h)$. There are two exceptional cases:  For $G(4,n,4)$, we take $G(p,q,h)=G(4,7,4)$ because $G(4,5,4)$ does not admits an embedding whose Euler genus is $\lceil f(4,5,4)\rceil$. For $G(m,n,0)$ with  $\tau_m=\tau_n=1$, we let $G(p,q,h)=G(3,3,0)$ since it must satisfy $q,q\geq 2.$
  A straightforward calculation yields the following choices for $G(p,q,h)$:
 \begin{align*}
 \intertext{When $k<m\leq n$, let}
G(p,q,h)&=\begin{cases}
    \displaystyle G(2, 2, 0), &\text{if}~ \tau_k=\tau_m=\tau_n=0,\\[0.3cm]
\displaystyle   G(2, 3, 0) ,& \text{if}~ \tau_k=0, \tau_m=0,\tau_n=1,\\[0.3cm]
\displaystyle   G(3, 4, 2) ,& \text{if}~ \tau_k=0, \tau_m=1,\tau_n=0,\\[0.3cm]
\displaystyle   G(3, 3, 2) ,& \text{if}~ \tau_k=0, \tau_m=1,\tau_n=1, ~\text{and}~k\geq 2,\\[0.3cm]
\displaystyle   G(3, 3, 0) ,& \text{if}~ \tau_k=0, \tau_m=1,\tau_n=1, ~\text{and}~k=0,\\[0.3cm]
\displaystyle   G(2, 2,1) ,& \text{if}~ \tau_k=1, \tau_m=0,\tau_n=0,\\[0.3cm]
\displaystyle   G(2, 3, 1) ,& \text{if}~ \tau_k=1, \tau_m=0,\tau_n=1,\\[0.3cm]
\displaystyle   G(3, 2, 1) ,& \text{if}~ \tau_k=1, \tau_m=1,\tau_n=0,\\[0.3cm]
\displaystyle G(3, 3, 1), &\text{if}~ \tau_k=\tau_m=\tau_n=1.
\end{cases}\\
\intertext{When $k=m<n$, let}
G(p,q,h)&=\begin{cases}
\displaystyle G(4,4,4), &\text{if}~ \tau_k=\tau_m=\tau_n=0,\\[0.3cm]
\displaystyle G(3,5,3), &\text{if}~ \tau_k=\tau_m=\tau_n=1,\\[0.3cm]
\displaystyle G(3,4,3) ,& \text{if}~ \tau_k=\tau_m=1, \tau_m=0,\\[0.3cm]
\displaystyle G(6,7,6) ,& \text{if}~ \tau_k=\tau_m=0, \tau_m=1, k=m\geq 6,\\[0.3cm]
\displaystyle G(4,7,4) ,& \text{if}~ \tau_k=\tau_m=0, \tau_m=1, k=m=4.
\end{cases}\\
\intertext{When $k=m=n$, let}
 G(p,q,h)&=\begin{cases}
\displaystyle G(4,4,4), &\text{if}~ \tau_k=\tau_m=\tau_n=0,\\[0.3cm]
\displaystyle G(5,5,5), &\text{if}~ \tau_k=\tau_m=\tau_n=1.
\end{cases}\\
\end{align*}

\begin{figure}[h]
\begin{center}
\setlength{\unitlength}{1.5mm}
  \begin{tikzpicture}
  \def\a{-1};
  \def\b{7};
 \coordinate [label=left:$b_7$]  (b7) at (\a+1.5,1.5);
\coordinate [label=left:$a_5$]  (a5) at (\a,4);
\coordinate [label=left:$b_6$]  (b6) at (\a,3);
\coordinate [label=left:$b_5$]  (b5) at (\a,2);
\coordinate [label=left:$a_6$]  (a6) at (\a,1);
\coordinate [label=right:$a_5$]  (a'5) at (\b,4);
\coordinate [label=right:$b_6$]  (b'6) at (\b,3);
\coordinate [label=right:$b_5$]  (b'5) at (\b,2);
\coordinate [label=right:$a_6$]  (a'6) at (\b,1);

\coordinate [label=left:$b_1$]  (b1) at (\a+3,4);
\coordinate [label=left:$a_2$]  (a2) at (\a+3,3);
\coordinate [label=left:$b_3$]  (b3) at (\a+3,2);
\coordinate [label=left:$a_4$]  (a4) at (\a+3,1);

\coordinate [label=right:$a_3$]  (a3) at (\b-3,4);
\coordinate [label=right:$b_4$]  (b4) at (\b-3,3);
\coordinate [label=right:$a_1$]  (a1) at (\b-3,2);
\coordinate [label=right:$b_2$]  (b2) at (\b-3,1);
\draw [thick]  (b1)--(a2) --(b3)--(a4)--({\a+3},0)(\a+3, 5)--(b1);
\draw [thick]  (\b-3,5)--(a3) (\b-3,0)--(b2)--(a1)--(b4)--(a3);
\draw [thick]  (b1)--(a3)(b3)--(a1)(a2)--(b4)(a4)--(b2);
\draw [thick] (a5)--(b1)(a6)--(\a+2,0)(\a+2,5)--(b1) (b5)--(a2)(b6)--(a2);
\draw [thick] (a'5)--(\b-2,5)(\b-2,0)--(b2)(a'6)--(b2)(b'5)--(a1) (b'6)--(a1);
\draw [thick](b7)--(a2)(b7)--(a4)(b7)--(a6);

\coordinate [label=center:{\small $\sim$}]  (H1) at (\a+1.5,2.8);

\coordinate [label=center:{\small $1$}]  (H2) at (\b-0.6,2.4);

\coordinate [label=center:{\small $1$}]  (H3) at (\a+4,4.5);

\coordinate [label=center:{\small $2$}]  (H4) at (\b-1,0.5);

\coordinate [label=center:{\small $2$}]  (H5) at (\a+4,2.5);

\filldraw [draw=black,fill=red!20]  (H1) circle (0.6cm);
\filldraw [draw=black,fill=blue!20]  (H2) circle (0.2cm);
\filldraw [draw=black,fill=blue!20]  (H3) circle (0.2cm);
\filldraw [draw=black,fill=blue!20]  (H4) circle (0.2cm);
\filldraw [draw=black,fill=blue!20]  (H5) circle (0.2cm);
\coordinate [label=center:{\small $1$}]  (H2) at (\b-0.6,2.4);
\coordinate [label=center:{\small $1$}]  (H3) at (\a+4,4.5);
\coordinate [label=center:{\small $2$}]  (H4) at (\b-1,0.5);
\coordinate [label=center:{\small $2$}]  (H5) at (\a+4,2.5);
\coordinate [label=center:{\small $\sim$}]  (H1) at (\a+1.5,2.8);

\draw [densely dotted,thick]  (b7)--({\a+1.5+0.6*cos(250)},{2.8+0.6*sin(250)});
\draw [densely dotted,thick]  (b7)--({\a+1.5+0.6*cos(270)},{2.8+0.6*sin(270)});
\draw [densely dotted,thick]  (b7)--({\a+1.5+0.6*cos(290)},{2.8+0.6*sin(290)});
\draw [densely dotted,thick]  (a5)--({\a+1.5+0.6*cos(90)},{2.8+0.6*sin(90)});
\draw [densely dotted,thick] (a3)--(\b,3.6) (\a,3.6)--({\a+1.5+0.6*cos(110)},{2.8+0.6*sin(110)});

\draw [densely dotted,thick] (a1)--(\b,3.4) (\a,3.4)--({\a+1.5+0.6*cos(130)},{2.8+0.6*sin(130)});

\draw [thick] (b'6)--({\b-0.6+0.2*cos(110)},{2.4+0.2*sin(110)});
\draw [thick] (b'6)--({\b-0.6+0.2*cos(50)},{2.4+0.2*sin(50)});
\draw [thick] (b'5)--({\b-0.6+0.2*cos(0)},{2.4+0.2*sin(0)});
\draw [thick] (b'5)--({\b-0.6+0.2*cos(270)},{2.4+0.2*sin(270)});

\draw [thick] (a3)--({\a+4+0.2*cos(0)},{4.5+0.2*sin(0)});
\draw [thick] (a3)--({\a+4+0.2*cos(270)},{4.5+0.2*sin(270)});
\draw [thick] (\a+3.8,5)--({\a+4+0.2*cos(130)},{4.5+0.2*sin(130)});
\draw [thick] (\a+4.4,5)--({\a+4+0.2*cos(70)},{4.5+0.2*sin(70)});
\draw [thick] (\a+3.8,0)--(a4);
\draw [thick] (\a+4.4,0)--(a4);

\draw [thick] (a'6)--({\b-1+0.2*cos(110)},{0.5+0.2*sin(110)});
\draw [thick] (a'6)--({\b-1+0.2*cos(20)},{0.5+0.2*sin(20)});

\draw [thick] (\b-0.8,0)--({\b-1+0.2*cos(-40)},{0.5+0.2*sin(-40)});
\draw [thick] (\b-1.2,0)--({\b-1+0.2*cos(250)},{0.5+0.2*sin(250)});
\draw [thick] (\b-0.8,5)--(a'5);
\draw [thick] (\b-1.2,5)--(a'5);

\draw [thick] (b3)--({\a+4+0.2*cos(180)},{2.5+0.2*sin(180)});
\draw [thick] (b3)--({\a+4+0.2*cos(270)},{2.5+0.2*sin(270)});
\draw [thick] (b4)--({\a+4+0.2*cos(90)},{2.5+0.2*sin(90)});
\draw [thick] (b4)--({\a+4+0.2*cos(0)},{2.5+0.2*sin(0)});

\draw [densely dashed, very thick,red]  (\a,0)--(\b,0)-- (\b,5)--(\a,5)--(\a,0);
\draw [densely dashed, very thick,red,->](\a, 2.2)--(\a, 2.4);
\draw [densely dashed, very thick,red,->](\b, 2.2)--(\b, 2.4);
\draw [densely dashed, very thick,red,->](2.4, 5)--(2.6, 5);
\draw [densely dashed, very thick,red,->](2.4, 0)--(2.6, 0);
 
\foreach \i in {1,2,3,4} {
\filldraw [draw=black,fill=red]  (a\i) circle (3pt);
\filldraw [draw=black,fill=blue]  (b\i) circle (3pt);
}
\foreach \i in {5,6} {
\filldraw [draw=black,fill=blue]  (b\i) circle (3pt);
\filldraw [draw=black,fill=blue]  (b'\i) circle (3pt);
\filldraw [draw=black,fill=red]  (a\i) circle (3pt);
\filldraw [draw=black,fill=red]  (a'\i) circle (3pt);
}
\filldraw [draw=black,fill=blue]  (b7) circle (3pt);
 
\end{tikzpicture}
\caption{{\small The embedding $\Pi_{6,7,6}$ with $\mathcal{E}(\Pi(6,7,6)=f(6,7,6)=7$. The red rectangle represents a Torus obtained by gluing the opposite sides of the rectangle along the directions indicated by the arrows; the circles with the same number are the ends of an handle.}
 }
\label{F11}
\end{center}
\end{figure}

Now, let $\Pi(p,q,h)$ denote an embedding of the graph $G(p,q,h)$ that satisfies $$\mathcal{E}(\Pi(p,q,h))=\lceil  f(p,q,h)\rceil.$$
Specifically, for
$(p,q,h)\in$ $\{(2,2,0)$,$(2,3,0)$,$(2,2,1)$,$(2,3,1)$,$(3,3,1)$,$(3,3,2)$,$(3,4,2)$, $(3,$ $4$,$3)$, $(3,5,3)$,$(4,4,4)\},$ 
$\Pi(p,q,h)$ is a planar embedding; $\Pi(3,3,0)$ is an embedding on the projective plane. 
These embeddings are straightforward to construct,  so we omit their illustrations. 
For $(p,q,h)\in$ $\{(5,5,5)$,$(4,7,4)$,$(6,7,6)\},$
the embeddings $\Pi(p,q,h)$ are shown in Figure \ref{F10} and  Figure \ref{F11}, respectively.

Therefore, by combing Theorem \ref{thm3}, the embedding 
$$\Pi(m,n,k)=\Pi(p,q,h)\oplus G(2,2,2)^a\oplus G(2,0,0)^b\oplus G(0,0,2)^c$$
is an embedding of $G(m,n,k)$ with 
$$\tilde{\gamma}(\Pi(m,n,k))=\lceil  f(m,n,k)\rceil.$$
  This completes the proof of Theorem \ref{thm2}.
 
 \section*{Acknowledgement} I would like to express my deepest gratitude to Professor Thomas Zaslavsky and Fengming Dong for polishing the English of this manuscript and for their insightful discussions.

\end{document}